\documentclass[10pt]{article}

\usepackage[top=1in,bottom=1.25in,left=1.25in,right=1.25in]{geometry}

\usepackage{graphicx}
\usepackage{ragged2e}
\usepackage{tabto}
\usepackage{amssymb}
\usepackage{float}
\usepackage{multirow}
\usepackage{amsmath}
\usepackage{tikz}
\usepackage{amsthm}
\usepackage{color}
\usepackage{url}
\usepackage[capitalise,noabbrev]{cleveref}
\usepackage{afterpage}
\usepackage{multicol}
\usepackage[font=footnotesize]{caption}
\usepackage[font=footnotesize]{subcaption}
\usepackage{marginnote}

\usepackage{lmodern}

\usepackage[sf,bf,medium]{titlesec}

\usepackage{bm}
\usepackage{dsfont}

\newcommand\mtx[1]{\bm{\mathsf{#1}}}
\newcommand\elem[1]{\mathsf{#1}}

\newcommand{\lp}{\left(}
\newcommand{\rp}{\right)}

\usepackage[sort,compress]{cite}

\newcommand\bbR{\mathbb R}

\newcommand\bx{\bm{x}}
\newcommand\bs{\bm{s}}
\newcommand\by{\bm{y}}

\newcommand\bc{\boldsymbol c}

\newcommand\In{\operatorname{inc}}

\newcommand\bn{\boldsymbol n}

\newcommand\bX{\boldsymbol X}

\newcommand\Npat{N_{\textrm{patches}}}
\newcommand\Sfar{\mathcal S_{\textrm{Far}}}
\newcommand\Snear{\mathcal S_{\textrm{Near}}}
\newcommand\Sself{\mathcal S_{\textrm{Self}}}

\newcommand\tlp{T_{\textnormal{LP}}}

\newcommand\tinit{T_{\textnormal{init}}}
\newcommand\tquad{T_{\textnormal{quad}}}

\newcommand\slp{S_{\textnormal{LP}}}

\newcommand\sinit{S_{\textnormal{init}}}
\newcommand\squad{S_{\textnormal{quad}}}

\newcommand\Nover{N_{\textnormal{over}}}
\newcommand\aavg{a_{\textnormal{avg}}}
\newcommand\amax{a_{\textnormal{max}}}

\newtheorem{remark}{\sffamily Remark}
\newtheorem{definition}{\sffamily Definition}

\newcommand{\cD}{\mathcal D}
\newcommand{\cS}{\mathcal S}
\newcommand{\cO}{\mathcal O}

\numberwithin{equation}{section}

\begin{document}


\begin{titlepage}

  \raggedleft
  {\texttt{Technical Report\\
    \today}}
  
  \hrulefill

  \vspace{4\baselineskip}

  \raggedright
  {\LARGE \sffamily\bfseries Fast multipole
  methods for the evaluation of layer potentials 
  with locally-corrected quadratures}
  

  \vspace{3\baselineskip}
 \vspace{\baselineskip}
 
  \normalsize Leslie Greengard\footnote{Research supported in part by 
  the Office of Naval Research under award
    number~\#N00014-18-1-2307.}\\
  \small \emph{Courant Institute, NYU \tabto{2in} Center for Computational
    Mathematics, Flatiron Institute\\
    New York, NY, 10012 \tabto{2in} New York, NY 10010}\\
  \texttt{greengard@cims.nyu.edu}

   \vspace{\baselineskip}
  \normalsize Michael O'Neil\footnote{Research supported in part by
    the Office of Naval Research under award
    numbers~\#N00014-17-1-2059,~\#N00014-17-1-2451, and~\#N00014-18-1-2307, and the Simons Foundation/SFARI (560651, AB).}\\
  \small \emph{Courant Institute, NYU\\
    New York, NY 10012}\\
  \texttt{oneil@cims.nyu.edu}
  
  \vspace{\baselineskip}
  \normalsize Manas Rachh\footnote{Corresponding author.}\\
  \small \emph{Center for Computational Mathematics, Flatiron Institute\\
    New York, NY 10010}\\
  \texttt{mrachh@flatironinstitute.org}
  \normalsize
  
  \vspace{\baselineskip}
  \normalsize Felipe Vico \footnote{Research supported in part by the Office of Naval Research under award number~\#N00014-18-2307, the Generalitat Valenciana under award number AICO/2019/018, and by the Spanish Ministry of Science and Innovation
(Ministerio Ciencia e Innovaci\'{o}n) under award number PID2019-107885GB-C32}\\
  \small \emph{Instituto de Telecomunicaciones y Aplicaciones Multimedia (ITEAM)\\ Universidad Polit`ecnica de Valencia\\
Valencia, Spain 46022}\\
  \texttt{felipe.vico@gmail.com}
  \normalsize

\end{titlepage}

\begin{abstract}
  While fast multipole methods (FMMs) are in widespread use for the rapid
  evaluation of potential fields governed by the Laplace, Helmholtz,
  Maxwell or Stokes equations, their coupling to high-order quadratures
  for evaluating layer potentials is still an area of active research. 
  In three dimensions, a number of issues need to be addressed, including
  the specification of the surface as the union of high-order patches,
  the incorporation of accurate quadrature rules for integrating singular 
  or weakly singular Green's functions on such patches, 
  and their coupling to the oct-tree data structures on which the FMM 
  separates near and far field interactions. Although the latter is 
  straightforward for point distributions, the near field for a patch 
  is determined by its physical dimensions, not the distribution of 
  discretization points on the surface.

  Here, we present a general framework for efficiently coupling
  locally corrected quadratures with FMMs, relying primarily on what
  are called generalized Gaussian quadratures rules, supplemented by
  adaptive integration.  The approach, however, is quite general and
  easily applicable to other schemes, such as Quadrature by Expansion
  (QBX).  We also introduce a number of accelerations to reduce the
  cost of quadrature generation itself, and present several numerical
  examples of acoustic scattering that demonstrate the accuracy,
  robustness, and computational efficiency of the scheme. 
  On a single core of an Intel i5 2.3GHz processor, a Fortran
  implementation of the scheme can generate near field quadrature 
  corrections for between 1000 and 10,000
  points per second, depending on the order of accuracy and the desired
  precision. A Fortran implementation of the algorithm described in
  this work is available at~\url{https://gitlab.com/fastalgorithms/fmm3dbie}.\\
  

\end{abstract}

\small

\tableofcontents

\normalsize
\newpage

\section{Introduction}

Over the last two decades, fast multipole methods (FMMs) 
and related hierarchical fast algorithms have become widespread
for computing $N$-body interactions in 
computational chemistry, astrophysics, acoustics, fluid dynamics, and 
electromagnetics. In the time-harmonic, acoustic setting, 
a typical computation of interest is the evaluation of
\begin{equation}
  F_m = \sum_{\substack{ n=1 \\ n\neq m}}^N
\sigma_n \, G_k(\bx_m,\bx_n) 
\label{fmmpt}
\end{equation}
where 
\begin{equation}
\label{eq:greenfunhelm}
G_k(\bx,\by) = \frac{e^{ik|\bx-\by|}}{4\pi|\bx-\by|}
\end{equation}
is the free-space Green's function for the Helmholtz equation 
\[ \Delta u + k^2 u = 0. \]
Direct calculation of~\cref{fmmpt} requires~$\cO(N^2)$ operations,
while the FMM requires~$\cO(N)$ work in the low frequency
regime~\cite{greengard_huang_etc} and~$\cO(N \log N)$ work in the high
frequency regime~\cite{wideband3d}.

When solving boundary value problems for partial differential
equations (PDEs) in three dimensions, such sums arise in the
discretization of layer potentials defined on a
surface~$\Gamma$. These potentials take the form:
\begin{equation}
  u(\bx) = \int_\Gamma K(\bx,\bx') \, \sigma(\bx') \, da(\bx').
\label{layerpotdef}
\end{equation}
In \eqref{layerpotdef}, $K(\bx,\bx')$ is a Green's function for the
PDE of interest, such as~\eqref{eq:greenfunhelm} or one of its
directional derivatives. As a result, the governing equation is
automatically satisfied, and it remains only to enforce the desired
boundary condition.  With a suitable choice for the kernel
$K(\bx,\bx')$, this often leads to a Fredholm integral equation of the
form
\begin{equation}
  \sigma(\bx) + \int_\Gamma K(\bx,\bx') \, \sigma(\bx') \, da(\bx') =
  f(\bx), \qquad \text{for } \bx \in \Gamma.
\label{fredholmeq}
\end{equation}
As we shall see below,
this can be discretized with high-order accuracy 
using a suitable Nystr\"om method \cite{Atkinson95,atkinson_1997}
\begin{equation}
  \sigma_i +  w_{ii} \sigma_i + 
\sum_{j\neq i} w_{ij} \, K(\bx_i,\bx_j) \, \sigma_j \,
  = f(\bx_i).
\end{equation}
Here, $\bx_i$ and $w_{ij}$ are the quadrature nodes and weights, 
respectively, while $\sigma_i$ is an approximation to the true
value~$\sigma(\bx_i)$. If the quadrature weights~$w_{ij}$ did not
depend on the target location, i.e.~$w_{ij} = w_j$, then the above sum
is a standard $N$-body calculation of the form~\eqref{fmmpt}.

Unfortunately, when the integral equation comes from 
a layer potential corresponding to an elliptic PDE,
the kernel~$K$ is typically singular or weakly singular 
so that simple high-order rules for smooth functions
fail.  Assuming the surface~$\Gamma$ is defined as the union of 
many smooth patches~$\Gamma_j$ (each with its own parameterization), 
high-order quadrature schemes require an analysis of the 
distance of the \emph{target} $\bx_i$ from each
{\em patch}.

More precisely, for 
a given target location~$\bx_i$ on the boundary, the integral
in \eqref{layerpotdef} or \eqref{fredholmeq}
can be split into three pieces: a self-interaction integral,
a near field integral, and a far field integral:
\begin{multline}
  \int_\Gamma K(\bx_i,\bx') \, \sigma(\bx') \, da(\bx') = 
  \int_{\text{Self}(\bx_i)} K(\bx_i,\bx' ) \, \sigma(\bx') \, da(\bx')
  \, + \\
  \int_{\text{Near}(\bx_i)} K(\bx_i,\bx' ) \, \sigma(\bx') \, da(\bx')
  \, +
  \int_{\text{Far}(\bx_i)} K(\bx_i,\bx' ) \, \sigma(\bx') \, da(\bx').
\end{multline}
This splitting into target-dependent regions is essential for
maintaining high-order accuracy. The region $\text{Self}(\bx_i)$ is simply
the patch on which $\bx_i$ itself lies. The integral over this patch
involves a singular integrand (due to the kernel~$K$). The \emph{Near} field
is defined precisely in Section~\ref{sec:local}, but consists of 
patches close enough to $\bx_i$ such that the 
integrand is \emph{nearly singular} even though it is formally
smooth.
The \emph{Far} region consists of all other patches, sufficiently far 
from~$\bx_i$ such that high-order quadratures  for smooth functions
can be applied. 

\begin{definition} \label{farfielddef}
Suppose that 
$\bx_i$ is in the far field of a patch $\Gamma_m$, and that 
\begin{equation}
\sum_{j=1}^M w_{j} \, K(\bx_i,\bs_j) \, \sigma_j \,  
  \approx
  \int_{\Gamma_m} K(\bx_i,\bx' ) \, \sigma(\bx') \, da(\bx')
\end{equation}
to the desired precision. Then 
$\bs_j$ and $w_{j}$ will be referred to as the {\em far field}
quadrature nodes and weights. Note that these are independent of $\bx_i$.
\end{definition}

A related task in the solution process is evaluating the
  computed solution $u(\bx)$ at target locations~$\bx$ which are
 off-surface, but possibly arbitrarily close to the surface. In
  this case as well, while the integrands are formally smooth,
  evaluating the integral presents a similar challenge owing to the
  nearly singular behavior of the integrand. The evaluation of the
  potential can still be split into two pieces in a similar manner: a
  near field integral and a far field integral
\begin{multline}
\int_\Gamma K(\bx,\bx') \, \sigma(\bx') \, da(\bx') = 
  \int_{\text{Near}(\bx)} K(\bx,\bx' ) \, \sigma(\bx') \, da(\bx') \\ +
  \int_{\text{Far}(\bx)} K(\bx,\bx' ) \, \sigma(\bx') \, da(\bx').
\end{multline}
The strategy outlined above can be used for evaluating the contributions of the
near and far regions for targets off-surface as well.

While methods exist for the Self and Near calculations,
the use of a fast algorithm such as the FMM 
requires coupling these somewhat complicated quadrature schemes to the 
Cartesian oct-tree data
structures that divide up space into a hierarchy of
regular, adaptively refined cubes. 
Unfortunately, the surface
patches~$\Gamma_m$ of the domain boundary~$\Gamma$
may be of vastly different sizes and do not, in general, conform to 
a spatial subdivision strategy based on the density of quadrature nodes
as points in $\mathbb{R}^3$. That is, many patches (curvilinear
triangles) $\Gamma_j$ are likely to
cross leaf node boundaries in the oct-tree
structure. 
If this were not the case, then far field interactions within the 
FMM could be handled by computing multipole expansions from entire 
patches, followed by direct calculation of the Self and Near
quadratures (analogous to the direct, near neighbor calculations in a
point-based FMM). 
Thus, one of the issues we address here concerns modifications of the FMM
so that speed is conserved for far field interactions, but in a manner 
where the overhead for near field quadrature corrections is modest, 
even when the surface patches are nonuniform.

Furthermore, while we will restrict our attention here to
Nystr\"om-style 
discretizations, the same concerns must be addressed for collocation and
Galerkin-type methods.  Similar issues arise when coupling adaptive
mesh refinement (AMR) data structures to complicated geometries using
Cartesian cut-cell methods~\cite{meshgenabm,
  Bell1994,johanssencolella}.



In the present paper, we develop an efficient algorithm
which allows for the straightforward coupling of adaptive FMM data 
structures with locally corrected quadrature schemes. 
Our goal is to achieve $\cO(N)$ or $\cO (N \log N)$ performance 
for surfaces with 
$O(N)$ discretizations points, depending on whether one is in the 
low or high frequency regime, respectively. Moreover,
we would like the constant implicit in this 
notation to be as close to the performance of point-based FMMs as
possible.
We will concentrate
on the use of generalized Gaussian quadrature
rules~\cite{bremer_2012c,bremer_2013,bremer-2015,bremer} for the self 
interactions on curvilinear triangles, and
adaptive integration for the Near region (nearly singular interactions). 
It is, perhaps, surprising that adaptive integration on surface patches
can be competitive with other schemes such as Quadrature By Expansion
(QBX), singularity subtraction, or coordinate
transformations 
\cite{bruno2001fast,bruno_garza_2020,malhotra19,erichsen1998quadrature,Siegel2018ALT,Wala2018,Wala2020,ying}.
The key is that we have developed a careful, precision and 
geometry-dependent hierarchy of 
interpolators on each patch, after which the adaptive integration step
is inexpensive when amortized over all relevant targets.
As a side-effect, our scheme also provides rapid access to 
entries of the fully discretized system matrix which is essential for
fast direct solvers.

\begin{remark}
Generalized Gaussian quadrature was already coupled to an FMM in \cite{bremer_2012c},
but the question of how to design a robust algorithm that is insensitive to large variation
in triangle dimensions was not directly addressed. In some sense, the present paper is devoted
to two separate issues raised in \cite{bremer_2012c}: the first is to accelerate adaptive
quadrature itself, and the second is to describe an FMM implementation that works for 
multi-scale discretizations.
\end{remark}

\begin{remark}
It is worth noting that most locally
corrected quadrature schemes, such as Duffy transformations \cite{Duffy},
are designed for a target that is mapped to
the origin of a local coordinate system or the vertex of a triangular
patch and, hence, are suitable only for self
interactions as defined above. Near interactions are not addressed. 
An exception is Quadrature by Expansion (QBX) which provides a
systematic, uniform procedure for computing layer potentials using
only smooth quadratures and extrapolation
\cite{klockner_2013,Siegel2018ALT,Wala2018}. These have been successfully
coupled to FMMs in \cite{Wala2018,Wala2020}.
Many aspects of the FMM modifications described here can be used in 
conjunction with QBX instead.
Another exception is 
Erichsen-Sauter rules \cite{erichsen1998quadrature}, which do include
schemes for adjacent panels, but appear to be best suited for 
modest accuracy.
\end{remark}

The paper is organized as follows: in Section~\ref{sec:prelim}, some
basic facts regarding polynomial approximation and integration on
triangles are presented. In Section~\ref{sec:surface}, we describe the
classical boundary integral equation for acoustic scattering from a
\emph{sound-soft} boundary, governed by the Helmholtz equation, as
well as discretization and integration methods for curvilinear
surfaces.  Section~\ref{sec:local} provides the algorithmic details
involved in locally corrected quadrature schemes.  The coupling of
these quadrature schemes to FMMs is presented in
Section~\ref{sec:fmmcoupling}, and numerical examples demonstrating the
performance of the scheme are presented in Section~\ref{sec:examples}.
Finally, in Section~\ref{sec:conclusions}, we discuss avenues for
further research, and the application of our scheme to fast direct
solvers.

\section{Interpolation and integration on
  triangles}
\label{sec:prelim}

For the sake of simplicity, we assume that we are given 
a surface triangulation represented as a collection of charts
$\bX^j = \bX^j(u,v)$, which map the standard right triangle
\begin{equation}
T_{0} = \{ (u,v): u\geq 0\,, v\geq 0 \,, u+v \leq 1 \}
\subset \bbR^{2} 
\label{eq:stdsimplex}
\end{equation}
to the surface patch $\Gamma_j$.
All discretization and integration is done over $T_{0}$, 
incorporating the mapping function $\bX^j$ and its derivatives as needed.

In this section, we summarize the basic polynomial
interpolation and integration
rules we will use for smooth functions $f: T_{0} \to \bbR$. 
A useful spectral basis is given by the
orthogonal polynomials on $T_{0}$, known as Koornwinder
polynomials~\cite{koornwinder_1975}. They are 
described  analytically by the formula:
\begin{equation}
  K_{nm}(u,v) = c_{nm} \, (1-v)^{m} \, P_{n-m}^{(0,2m+1)}(1-2v)  \,
  P_{m}\lp \frac{2u+v-1}{1-v} \rp , \qquad m\leq n ,
\end{equation}
where~$P_{n}^{(a,b)}$ is the Jacobi polynomial of degree~$n$ with
parameters~$(a,b)$, $P_{m}$~is the Legendre polynomial of degree~$m$,
and~$c_{nm}$ is a normalization constant such that
\begin{equation}
  \int_{T_{0}} \vert K_{nm}(u,v) \vert^2 \, du \, dv = 1.
\end{equation}
For convenience, our definition is slightly different from 
that in~\cite{koornwinder_1975}.

It is easy to see that
there are~$n_{p} = p(p+1)/2$ Koornwinder polynomials of total degree less
than~$p$. By a straightforward change of variables, and
using the orthogonality relationships for Legendre and Jacobi
polynomials, it is easy to show that
\begin{equation}
  \int_{T_0} K_{nm}(u,v) \, K_{n'm'}(u,v) \, du \, dv = 0, \qquad \text{for }
  n\neq n' \text{ and }  m\neq m'.
\end{equation}
The Koornwinder polynomials form a complete basis for~$L^2(T_0)$, and
can easily be used to approximate smooth functions on~$T_0$ to
arbitrary precision.

\subsection{Polynomial approximation and integration}
\label{sec:approx}

As in standard spectral approximation methods 
for functions defined on intervals or tensor products of intervals
\cite{trefethen2013approx},
smooth functions~$f$ defined on~$T_0$ can be interpolated,
approximated, and integrated using a Koornwinder polynomial basis.

To this end, suppose that~$f$ is defined by a $p$th-order Koornwinder 
expansion with coefficients $c_{nm}$:
\begin{equation}
f(u,v) = \sum_{n=0}^{p-1} \sum_{m=0}^{n} c_{nm} \, K_{nm}(u,v).
\end{equation}
Then, the square matrix~$\mtx{U}$ that maps the~$n_p$ coefficients in
the above expansion to values of~$f$ at a selection of~$n_p$
interpolation nodes, denoted by~$(u_j,v_j) \subseteq T_0$, has elements
\begin{equation}
  \elem{U}_{nm,j} = K_{nm}(u_j,v_j) \label{eq:defumat}.
\end{equation}
Let the matrix~$\mtx{V} = \mtx{U}^{-1}$. Then~$\mtx{V}$ maps values
of~$f$ at the interpolation nodes~$(u_j,v_j)$ to coefficients in a
Koornwinder polynomial expansion.

Suppose now that $f:T_{0}\to \bbR$ is an arbitrary smooth function
(not necessarily a polynomial) and let the values of~$f$ at the
interpolation points be denoted by~$f_{i} = f(u_{i},v_{i})$. 
Then, a $p$th-order approximation to $f$ is given by
\begin{equation}
    f(u,v) \approx \sum_{n=0}^{p-1}\sum_{m=0}^{n} c_{nm} \, K_{nm}(u,v)  ,
\end{equation}
where
\begin{equation}
  \label{eq:v2c}
    c_{nm} = \sum_{j=1}^{n_{p}} \elem{V}_{(nm),j} \, f_{j}  .
\end{equation}
We define the \emph{conditioning of the interpolation procedure} as
the condition number of the matrices~$\mtx{U}$ or $\mtx{V}$. Much like
interpolation operators on the interval, these matrices are not
well-conditioned for arbitrary selections of nodes~$(u_j,v_j)$, as we
will briefly discuss in the next two sections.

An $n$-point quadrature rule for computing the integral of a
function~$f$ on~$T_0$ is a collection of nodes and weights, 
$(u_{i},v_{i})$, $w_{i}$, $i=1,2,\ldots n$, such that
\begin{equation}
  \begin{aligned}
    \int_{T_{0}} f(u,v) \, du \, dv &= \int_{0}^{1} \int_{0}^{1-u}
    f(u,v) \, du \,dv \\
    &\approx \sum_{i=1}^{n}  w_{i} \, f_{i},
  \end{aligned}
\end{equation}
where $f_{i} = f(u_{i},v_{i})$ is the value of $f$ at
the $i$th quadrature node. The accuracy of such a quadrature rule is very
dependent on the choice of the quadrature nodes; if the rule is to be
exact for a selection of~$n$ functions, then the values of~$w_i$ are
determined wholly by the selection of~$(u_i,v_i)$.
Since the node selection provides additional degrees of freedom,
Gaussian-type quadrature rules are possible.
On the triangle, a quadrature rule would be perfectly Gaussian if 
it integrated~$3n$ functions
exactly, since there are~$3n$ parameters (the two
coordinates of the nodes, and the weights). In one dimension, it is 
well-known that choosing the nodes as the roots of a suitable orthogonal
polynomial leads to a perfect Gaussian rule \cite{trefethen2013approx}.
In two dimensions, such
perfect rules do not exist, but 
approximately Gaussian quadrature rules can be constructed.

\subsection{High-order quadrature rules on the simplex~$T_0$}

First described in~2010~\cite{xiao2010numerical}, what we will
refer to as Xiao-Gimbutas quadratures are a set of Gaussian-like
rules obtained through the solution of a nonlinear least-squares problem
using Newton's method. The
resulting nodes are contained in the interior of~$T_0$, and all the
weights are positive. Various kinds of symmetry can also be specified.
For a given $p>0$,
these rules are designed to use the minimum
number of nodes with positive weights so that the
resulting quadrature rule is \emph{exact} for all polynomials of
total degree~$<p$. As noted above, there are~$n_p = p(p+1)/2$
such polynomials. While not perfect Gaussian rules,
the Xiao-Gimbutas quadratures achieve remarkably high-order.
A rule with 48 weights and nodes, for example, is exact for polynomials
of degree 16, of which there are 
$n_{16} = 136$. The rule has only $3 \times 48 = 144$ free parameters.

For the sake of convenience we would like the quadrature nodes to
serve as interpolation/approximation nodes as well. There are,
however, far fewer Xiao-Gimbutas nodes than functions we would like to
interpolate (namely $n_p$).  Instead of using an even higher order
Xiao-Gimbutas rule, with at least $n_p$ nodes, we choose an
alternative quadrature scheme, introduced by Vioreanu and Rokhlin in
2014~\cite{vioreanu_2014}.  The Vioreanu-Rokhlin nodes of order~$p$,
obtained via a similar optimization procedure, are a collection
of~$n_p$ nodes which can be used simultaneously for high-order
polynomial interpolation, approximation, and integration on
$T_{0}$. We refer the reader to~\cite{vioreanu_2014} for a thorough
discussion. For our purposes, it suffices to note that the
interpolation operators computed using these nodes are extremely
well-conditioned.

As a quadrature rule, the nodes and weights are Gaussian-like; they
have positive weights and integrate more polynomials than there are
nodes in the quadrature. For example, the Vioreanu-Rokhlin rule that
interpolates polynomials of degree $p=16$, with $n_{16}=136$ nodes,
integrates all $378$ polynomials of degree up to $p' = 27$.  A perfect
Gaussian rule would integrate~$3n_p = 408$ functions exactly.  The
relationship between the order of the interpolation scheme~$p$ and the
order of the quadrature $p'$ is somewhat complicated, and obtained
empirically~\cite{vioreanu_2014}.

\section{Acoustic scattering from a sound-soft boundary}
\label{sec:surface}

Let~$\Omega$ be a bounded region in~$\bbR^{3}$, with smooth
boundary~$\partial\Omega = \Gamma$.  Given a function~$f$
defined on~$\Gamma$, 
a function~$u$ defined in $\bbR^{3} \setminus \Omega$ is said to 
satisfy the exterior Dirichlet problem for the Helmholtz
equation if
\begin{equation}\label{eq:extdir}
  \begin{aligned}
    (\Delta + k^2) \, u &=  0 &\qquad &\text{in } \bbR^{3}
    \setminus \Omega,\\
    u &= f & &\text{on } \Gamma, \\
    \lim_{r \to \infty} \, r \left(\frac{\partial u}{\partial r} - ik u
    \right) &= 0. & &
\end{aligned}
\end{equation}
In acoustics, Dirichlet problems such as this arise
when~$\partial\Omega$ is \emph{sound-soft} and $f = - u^{in}$, where
$u^{in}$ is an impinging acoustic wave. A standard approach for
solving the Dirichlet problem is to let
\begin{equation}
\label{eq:urep}
u = \cD_{k}[\sigma] - ik \, \cS_{k}[\sigma],
\end{equation}
where~$\sigma$ is an unknown density function defined
on~$\Gamma$.  Here~$\cS$ and~$\cD$ are the single layer and double
layer operators, respectively, given by
\begin{align}
  \cS_{k}[\sigma](\bx)
  &= \int_{\Gamma} G_{k}(\bx,\by) \, \sigma(\by) \, da(\by),  \label{eq:sldef} \\
  \cD_{k}[\sigma](\bx)
  &= \int_{\Gamma}  \lp \bn(\by) \cdot \nabla_{\by}G_{k}(\bx,\by) \rp \,  
    \sigma(\by) \, da(\by), \label{eq:dldef}
\end{align}
where $\bn(\by)$ is the outward normal at~$\by \in \Gamma$, and
$G_k(\bx,\by)$ is given by~\eqref{eq:greenfunhelm}.
The representation~\eqref{eq:urep} automatically satisfies the
Helmholtz equation and the radiation condition in~\eqref{eq:extdir}.
Imposing the boundary condition, and using standard jump relations for
layer potentials~\cite{colton_kress}, we obtain the following
second-kind integral
equation along~$\Gamma$ for the density $\sigma$:
\begin{equation}
\label{eq:inteq}
\frac{1}{2}\sigma(\bx) + \cD_{k}[\sigma](\bx) \textendash ik \cS_{k}[\sigma](\bx) 
= f(\bx), \qquad \bx \in \Gamma.
\end{equation}
This involves a slight abuse of notation: for $\bx \in \Gamma$,
$\cD_{k}[\sigma](\bx)$ should be
evaluated in the principal value sense.

When solving~\eqref{eq:inteq}, the accurate evaluation of the layer
potentials $\cS_{k}[\sigma]$, $\cD_{k}[\sigma]$ on~$\Gamma$ is essential
for either direct or iterative solvers. We will focus here
on the evaluation of the single layer
potential $\cS[\sigma]$, assuming~$\sigma$ is known.
Only minor modifications are needed
to address the double layer
potential, as well as other scalar or vector-valued layer
potentials that arise in electrostatics, elastostatics, viscous
flow, or electromagnetics.

\subsection{Surface parameterizations}

While some simple boundaries (such as spheres, ellipsoids and tori)
can be described by global parameterizations, in general it is
necessary to describe a complicated surface~$\Gamma$ as a collection
of surface patches, each of which is referred to as a
\emph{chart}. The collection of charts whose union defines $\Gamma$
will be referred to as an \emph{atlas}.

More precisely, we assume that 
the surface is the disjoint union of patches~$\Gamma_j$
\begin{equation}
  \Gamma = \cup_{j=1}^{\Npat} \, \Gamma_{j},
\end{equation}
and that the patch~$\Gamma_{j}$ is parameterized by a non-degenerate
chart~$\bX^j: T_{0} \to \Gamma_{j}$, where $T_{0}$ is the standard
simplex~\eqref{eq:stdsimplex}. Given these charts~$\bX^j$,
a local coordinate system can be defined on patch~$\Gamma_j$ by
taking its partial derivatives. For this, we define
\begin{equation}
  \bX^j_u \equiv \frac{\partial \bX^j}{\partial u}, \qquad
  \bX^j_v \equiv \frac{\partial \bX^j}{\partial v}, \qquad
  \bn^j \equiv \bX^j_u \times \bX^j_v.
\end{equation}
Finally, we assume that the triplet~$\bX^j_u$, $\bX^j_v$, $\bn^j$ forms
a right-handed coordinate system with~$\bn^j$ pointing into the
unbounded region $\bbR^3 \setminus \Omega$. 
In general, these vectors are neither
orthogonal nor of unit length. The area element on the
patch~$\Gamma_j$ is determined by the Jacobian~$J^j$,
\begin{equation}\label{eq:loc}
  \begin{aligned}
    da(\bX^j) &= \vert \bX^j_u \times \bX^j_v \vert \, du \, dv \\
    &= J^j(u,v) \, du \, dv.
  \end{aligned}
\end{equation}

\subsection{Discretization and integration}
\label{sec:disc}

If~$f$ is a function defined on~$\Gamma$, then
its integral can be decomposed as a sum over patches:
\begin{equation}\label{eq:intsplit}
  \begin{aligned}
    \int_{\Gamma} f(\bx) \, da(\bx)
    &= \sum_{j=1}^{\Npat} \int_{\Gamma_{j}} f(\bx) \, da(\bx) \\
    &= \sum_{j=1}^{\Npat} \int_{T_{0}} f\left(\bX^{j}(u,v)\right) \,
    J^{j}(u,v) \, du  dv \\
    &= \sum_{j=1}^{\Npat} \int_{u=0}^{1} \int_{v=0}^{1-u}
        f\left(\bX^{j}(u,v)\right) \, J^{j}(u,v) \,  du dv,
\end{aligned}
\end{equation}
where the Jacobian is given in~\eqref{eq:loc}.

If, in addition, $f$ is smooth, then each of the integrals on~$T_{0}$
in~\eqref{eq:intsplit} can be approximated using Xiao-Gimbutas or
Vioreanu-Rokhlin quadrature rules, as discussed in
Section~\ref{sec:prelim}. Using the latter. we have
\begin{equation}
    \int_{\Gamma} f(\bx) \, da(\bx) \approx
\sum_{j=1}^{\Npat} \sum_{\ell=1}^{n_p} w_\ell \, 
f\left(\bX^{j}(u_\ell,v_\ell)\right) \, J^{j}(u_\ell,v_\ell),
\end{equation}
where~$n_p$ is the number of nodes in the quadrature, which varies
depending on the desired order of accuracy.

For a patch~$\Gamma_{j}$ and a target $\bx$, however, the integrand
appearing in $S_{k}[\sigma](\bx)$ is only smooth when $\bx$ is in the
far field. When $\bx$ is either on $\Gamma_j$ or nearby, we will need
to modify our quadrature approach, as described at the outset.
Furthermore, in practice, we will only be given approximations of the
charts~$\bX^j$ and the function~$f$ (or the density~$\sigma$)
on each patch to finite order.

In what follows, we define a $p$th-order approximation as one for
which the truncation error is~$\mathcal O(h^p)$, where $h$ is, say,
the diameter of the patch $\Gamma_j$.  For a scalar function~$f$, it
can be approximated to $p$th-order in~$L^2(T_0)$ using Koornwinder
polynomials as
\begin{equation}
  f(u,v) \approx \sum_{n+m < p} f_{nm} \, K_{nm}(u,v).
\end{equation}
The coefficients~$f_{nm}$ can be computed using the
values-to-coefficients matrix~$\mtx{V}$, as described in
Section~\ref{sec:approx}.
The charts~$\bX^j$ will generally be approximated using a
vector version of the above formula:
\begin{equation}
  \begin{aligned}
  \bX^j(u,v) &\approx \sum_{n+m <p}
  \begin{pmatrix}
    x^j_{nm} \\
    y^j_{nm} \\
    z^j_{nm}
  \end{pmatrix} K_{nm}(u,v) \\
  &=\sum_{n+m <p} \bx_{nm}^j \, K_{nm}(u,v)
  \end{aligned}
\end{equation}
It is important to note that even if the charts~$\bX^j$ are
approximated to accuracy~$\epsilon$, i.e.
\begin{equation}
  \bX^j(u,v) = \sum_{n+m <p} \bx_{nm}^j \, K_{nm}(u,v) + \mathcal O(\epsilon),
\end{equation}
it does not follow that the Jacobian~$J^j$ will be evaluated to
precision~$\epsilon$ as well. The function~$J^j$ is non-linear and usually
requires a higher-order approximation than the individual components 
of the chart itself. This cannot be avoided unless analytic derivative 
information is provided for each patch~$\Gamma_j$.
This affects the accuracy of numerical approximations to surface integrals
for a fixed set of patches (but not the asymptotic convergence rate).

\begin{remark}
  In the following, we use the same order of discretization for
  representing the layer potential densities $\sigma$ and the surface
  information, i.e. the charts $\bX^{j}$ and their derivatives
  $\bX^{j}_{u}$, and $\bX^{j}_{v}$. While this choice is made for
  convenience of notation and software implementation, our approach to
  evaluating layer potentials extends in a straightforward manner to
  the case where different orders of discretization are used for
  representing the surface information and layer
  potentials densities.
\end{remark}

\section{Locally corrected quadratures}
\label{sec:local}

For a target location $\bx \in \Gamma_j$, let us first consider the 
{\em self interaction}
\begin{equation}
\Sself[\sigma](\bx) =
   \int_{\Gamma_{j}} G_k(\bx,\by) \, \sigma(\by) \, da(\by).
\label{sselfdef}
\end{equation}
In~\cite{bremer_2012c,bremer_2013}, the authors designed quadrature
rules for exactly this purpose, under the assumption that~$\sigma$
and~$J^j$ are well-approximated by polynomials (and therefore
representable by Koornwinder expansions).  The quadrature schemes in
these papers involve a rather intricate set of transformations but
yield a set of precomputed tables which, when composed with the chart
$\bX^j$, yield the desired high-order accuracy. Briefly, for any
  $\bx \in \Gamma_{j}$, and all $\sigma$ of the form
  $\sigma(u,v) = \sum_{nm} c_{nm} K_{nm}(u,v)$, there exist $N(\bx)$
  nodes $(u_{\bx,\ell},v_{\bx,\ell}) \in T_{0}$ and associated
  quadrature weights $w_{\bx,\ell}$, such that
\begin{multline}
\label{eq:errestself}
\bigg | \int_{u=0}^{1} \int_{v=0}^{1-u} G_{k}(\bx, \bX^{j}(u,v)) \,
  \sigma(u,v)
  \, J^{j}(u,v) \, du \, dv   \\
 - \sum_{\ell=1}^{N(\bx)} G_{k}(\bx,
  \bX^{j}(u_{\bx,\ell},v_{\bx,\ell})) \,
  \sigma(u_{\bx,\ell},v_{\bx,\ell}) \,
  J^{j}(u_{\bx,\ell},v_{\bx,\ell}) w_{\bx,\ell} \bigg | \leq
\varepsilon \cdot \| \sigma \|_{\mathbb{L}^{2}(\Gamma_{j})}.
\end{multline}

As mentioned in the introduction, in the original
paper~\cite{bremer_2012c}, which was focused on quadrature design, a
simple coupling to FMMs was mentioned that relied on the underlying
discretization being uniformly high-order. Near field interactions
were done using \emph{on-the-fly} adaptive integration. In their
subsequent paper~\cite{bremer_2013}, this type of expensive adaptive
integration was used for all non-self interactions (i.e. no FMM-type
acceleration was used at all).  Such an approach cannot be directly
accelerated with standard FMMs since the effective quadrature weights
are functions of both the source and target locations.

Recall that
for $\bx \in \bbR \setminus \Gamma$, we split the single layer potential
$\cS_{k}[\sigma](\bx)$ into two pieces:
\begin{equation}
\label{eq:near-far-split-off}
  \begin{aligned}
    \cS_{k}[\sigma](\bx) &= \int_{\Gamma} G_k(\bx,\by) \, \sigma(\by)
    \, da(\by) \\
    &=\sum_{\ell=1}^{\Npat} \int_{\Gamma_{\ell}} G_k(\bx,\by) \,
    \sigma(\by) \,  da(\by)  = 
    \Snear[\sigma](\bx) + \Sfar[\sigma](\bx) \, ,
\end{aligned}
\end{equation}
and when $\bx$ lies {\em on} the boundary, say on patch $\Gamma_j$, then
the single layer potential
$\cS[\sigma](\bx)$ is split into three pieces:
\begin{equation}
\label{eq:near-far-split-on}
    \cS_{k}[\sigma](\bx) = 
    \Sself[\sigma](\bx) + \Snear[\sigma](\bx) + \Sfar[\sigma](\bx),
\end{equation}
where $\Sself[\sigma](\bx)$ is defined in 
\eqref{sselfdef}, and the near and far regions associated with a target and the corresponding definitions 
of $\Snear[\sigma](\bx)$ and $\Sfar[\sigma](\bx)$ are described below.

\begin{figure}[t]
  \centering
  \includegraphics[width=0.5\linewidth]{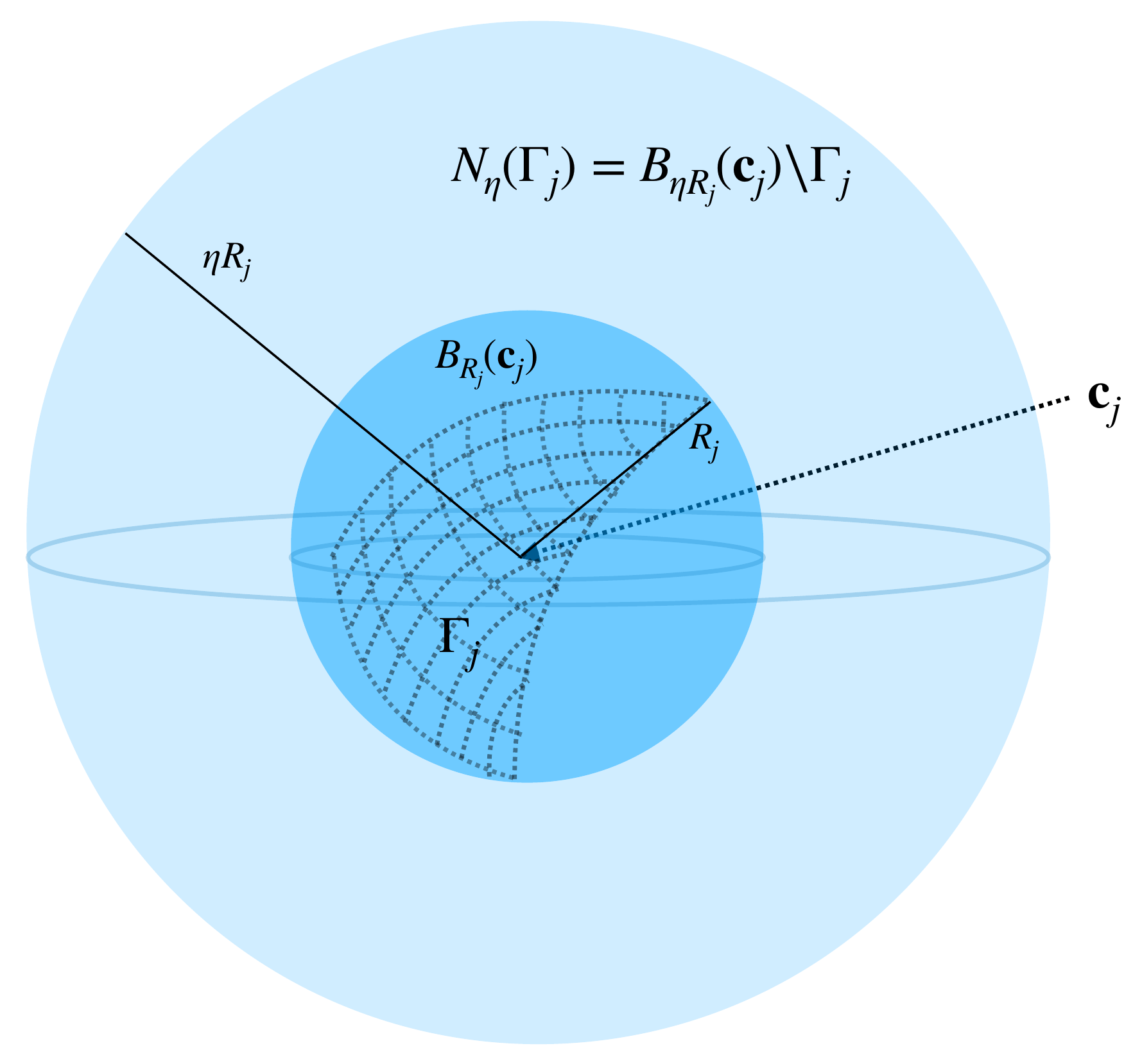}
  \caption{The smallest sphere containing surface patch~$\Gamma_j$
    centered at~${\bf c}_j$ and the near field region
    $N_\eta(\Gamma_j)$.  $\eta>1$ is a free parameter whose selection
    is based on the order of accuracy of the far field quadrature.}
  \label{fig:near-field-def}
\end{figure}

For this, it turns out to be easier to first take the point of view of a 
patch rather than a target.
Let~$\bc_{j}$ denote the
centroid of the patch~$\Gamma_{j}$,
\begin{equation}
  \begin{aligned}
    \bc_{j} &= \int_{\Gamma_{j}} \bx \, da(\bx) \\
    &= \int_{T_{0}} \bX^{j}(u,v) \, du \, dv ,
  \end{aligned}
\end{equation}
and let 
\begin{equation}
  \label{eq:rjdef}
  R_{j} = \min_{R>0} \{ R \mid \, \Gamma_{j} \subset B_{R}(\bc_{j}) \}, 
\end{equation}
where~$B_{R}(\bc_j)$ is the ball of radius~$R$ centered at~$\bc_j$.
That is, $B_{R_j}(\bc_j)$ is the ball of minimal radius containing the
patch $\Gamma_j$. 
Letting $\eta > 1$ be a free parameter for the moment, we define the 
near field of the patch $\Gamma_{j}$, denoted by 
$N_{\eta}(\Gamma_{j})$, to be the set of points that do not lie on
$\Gamma_j$ but are within the ball
$B_{\eta R_j}(\bc_j)$ (see Fig.~\ref{fig:near-field-def}). Thus,
\begin{equation}
N_{\eta}(\Gamma_{j})
= \{\bx \in \bbR^3 \setminus \Gamma_j \mid \, 
d(\bc_{j},\bx) \leq \eta R_{j} \} \, .
\end{equation}
Given the collection of near field regions of the form
$N_{\eta}(\Gamma_{j})$, let $T_{\eta}(\bx)$ denote the dual list:
that is, 
the collection of patches $\Gamma_{j}$ for which the point
$\bx \in N_{\eta}(\Gamma_{j})$,
\begin{equation}
    T_{\eta}(\bx) = \{ \Gamma_{j} \mid \bx \in N_{\eta}(\Gamma_{j}) \}.
\end{equation}
Similarly, for $\bx \in \bbR^3 \setminus \Gamma$, we denote the far 
field of $\bx$ by
\begin{equation}
    F_{\eta}(\bx) = \{ \Gamma_{j} \mid  \bx \notin N_{\eta}(\Gamma_{j}) \}.
\end{equation}
When $\bx \in \Gamma_i$ is a boundary point, we let
\begin{equation}
    F_{\eta}(\bx) = \{ \Gamma_{j}, j \neq i \mid  \bx \notin N_{\eta}(\Gamma_{j}) \}.
\end{equation}

Then referring to the near-far split of the layer potential for targets off and on-surface described in~\cref{eq:near-far-split-off}, 
and~\cref{eq:near-far-split-on}, respectively, the near and far part of the layer potentials $\Snear$, and $\Sfar$ are given by
\begin{equation}
\Snear[\sigma](\bx) =
    \sum_{\Gamma_{\ell} \in T_{\eta}(\bx)}
    \int_{\Gamma_{\ell}} G_k(\bx, \by) \,  \sigma(\by) \, da(\by) 
\end{equation}
and
\begin{equation}
\Sfar[\sigma](\bx) = 
    \sum_{\Gamma_{\ell} \in F_{\eta}(\bx)}
    \int_{\Gamma_{\ell}}
    G_k(\bx, \by) \, \sigma(\by) \, da(\by).
\end{equation}

As noted in the beginning of the section, $\Sself[\sigma](\bx)$ can be
computed using the generalized Gaussian quadratures of 
\cite{bremer_2012c,bremer_2013}.
By virtue of their separation from the source patches,
all of the integrands in $\Sfar$ are smooth 
and can be computed using
either Vioreanu-Rokhlin or Xiao-Gimbutas quadrature rules,
with weights that are independent of 
the target location~$\bx$. The accuracy of these rules, which is
affected by the free parameter $\eta$, is discussed
in Section~\ref{sec:etachoice}.

It remains only to develop an efficient scheme for evaluating the
integrals which define~$\Snear[\sigma]$.  At present, there do not exist
quadrature rules that are capable of simultaneously accounting for the
singularity in the Green's function and the local geometric variation
in an efficient manner.  The approach developed below involves a
judicious combination of precomputation and a greedy, adaptive
algorithm applied, for every target point~$\bx$, to
each~$\Gamma_{\ell}\in T_{\eta}(\bx)$.  Once the near field
quadratures have been computed, they can be saved using only $O(N)$
storage.  When solving an integral equation iteratively, this can be
used to accelerate subsequent applications of the integral operator.

\begin{remark}
The evaluation of near field quadratures does {\em not} affect the
overall complexity of computing layer potentials, assuming~$\eta$ is
not too large. This follows from the fact that, for each target~$\bx$,
there are only~$O(1)$ patches contained in~$T_{\eta}(\bx)$.  Thus, the
cost of all near field contributions is of the order~$O(N)$.  Since
there can be several patches in the near field, however, this
computation tends to be the rate limiting step in the overall
quadrature generation procedure.
\end{remark}

\begin{remark}
Without entering into a detailed literature review, it should be noted
that coordinate transformation methods such as those in 
\cite{bruno2001fast,malhotra19,erichsen1998quadrature,ying} can also be used
for computing near field interactions for surface targets. However,
these methods don't apply easily to off-surface evaluation.
An alternative to our procedure is
Quadrature by Expansion (QBX)
\cite{klockner_2013,Siegel2018ALT,Wala2018,Wala2020},
which handles singular and nearly-singular
integrals in a unified manner and (like the method of this paper) 
works both on and off surface.
There are distinct trade-offs to be made in QBX-based schemes and the
scheme presented here. In the end, the 
best method will be determined by accuracy, efficiency and ease of use.
At present, we have found the adaptive quadrature approach
to be the most robust and fastest in terms of overall performance. 
\end{remark}


\subsection{Selecting the near field cutoff \label{sec:etachoice}}

In this section, we discuss the choice of the parameter $\eta$, which
defines the near field for each patch~(see
Fig.~\ref{fig:near-field-def}).  Once $\eta$ is fixed, accuracy
considerations will determine whether the interpolation nodes used on
each patch are sufficient for accurate calculation of the far field
$\Sfar$, or whether we will need to increase the order of the far
field quadrature.

As $\eta$ increases, the near field for each patch obviously grows, so
that the number of targets for which we will apply specialized
quadrature increases. Since we would like to store these near field
quadratures for the purpose of repeated application of the integral
operator, both the storage and CPU requirements also grow
accordingly. On the other hand, as $\eta$ decreases, the integrand in
$\Sfar$ becomes less smooth, and the number of quadrature nodes needed
to achieve the desired precision in the far field will grow.  If it
exceeds the number of original nodes $n_p$ to achieve $p$th-order
convergence, we will have to {\em oversample} the layer potential
density~$\sigma$.  That is, we will have to interpolate~$\sigma$ to a
larger number of quadrature nodes.  This increases the computational
cost of evaluating the far field via the FMM.  Balancing the cost of
the near field and far field interactions sets the optimal value for~$\eta$.

Based on extensive numerical experiments, we have found that for the
highest order methods ($p>8$), $\eta = 1.25$ works well.  For orders
of accuracy $4<p\leq 8$, we recommend $\eta = 2$, and for the lowest
orders of accuracy $p\leq 4$, we recommend $\eta = 2.75$.  

\subsection{Oversampling via
    $p-$refinement} \label{sec:oversamp} For a patch
  $\Gamma_{j}$, suppose that an order $q$ Vioreanu-Rokhlin quadrature
  rule is the smallest order quadrature rule which accurately computes
  the 
  contribution of $\Gamma_{j}$ to all targets
  $\bx \not \in N_{\eta}(\Gamma_{j})$ to the desired precision
  $\varepsilon$. Then, the oversampling factor for $\Gamma_{j}$ is
  defined as the ratio of $q(q+1)/(p(p+1))$.  For a given $\eta$,
rather than trying to estimate the oversampling factor needed for
$\Sfar$ analytically, we compute the far field quadrature order $q$
needed for a specified precision numerically.  For each patch
$\Gamma_{j}$, we first identify the $10$ farthest targets in
$N_{\eta}(\Gamma_{j})$.  If the cardinality of $N_{\eta}(\Gamma_{j})$
is less than $20$, we choose the farthest $|N_{\eta}(\Gamma_{j})|/2$
targets from the list, and append $15$ randomly chosen targets on the
boundary of the sphere $\partial B_{\eta R_{j}} (\boldsymbol{c}_{j})$.
We denote this set of targets by $F(\Gamma_{j})$.  For the $nm$-th
Koornwinder polynomial $K_{nm}$, let
\begin{equation}
  I^j_{nm} (\bx) = \int_{T_{0}} G_{k}(\bx, \bX^{j}(u,v)) \,
  K_{nm}(u,v) \, J^{j}(u,v) \, du \, dv ,
\end{equation}
and let $\tilde{I}^j_{nm,q}(\bx)$ 
denote the approximation to the integral computed 
using the $q$th-order Vioreanu-Rokhlin quadrature. 
Then, the far field order $q_{j}$ for patch $\Gamma_{j}$ is 
chosen according to the following criterion: $q_j$ is the 
smallest $q$ such that all of the integrals 
$\tilde{I}^j_{nm,q}(\bx)$, for $0\leq m\leq n \leq p$
and $\bx \in F(\Gamma_{j})$, agree to a prescribed
tolerance~$\varepsilon d_{j}/\| \mtx{V} \|$ with the corresponding integrals 
obtained from using a~$(q+1)$th-order Vioreanu-Rokhlin quadrature. Here, $d_{j}$ is given by 
\begin{equation}
d_{j} = \min_{\ell=1\ldots n_{p}} \sqrt{J^{j}(u_{\ell},v_{\ell}) w_{\ell}},
\end{equation}
where $u_{\ell},v_{\ell}$ are the order $p$ Vioreanu-Rokhlin nodes and $w_{\ell}$ are the corresponding weights and $\| \mtx{V} \|$ is the operator norm of the values to interpolation matrix $\mtx{V}$ defined in~\cref{sec:approx}. That
is to say,
\begin{equation}
  \label{eq:qj}
q_{j} = \min_{q} \text{ such that } \max_{\bx \in F(\Gamma_{j})} \sqrt{\sum_{n+m<p} |\tilde{I}^j_{nm,q}(\bx) - \tilde{I}^j_{nm,q+1}(\bx)|^2} \leq \varepsilon \frac{d_{j}}{\| \mtx{V} \|} \,.
\end{equation}

This seemingly arbitrary choice of scaling $\varepsilon$ allows us to obtain an estimate for the relative error of the contribution of $\Gamma_{j}$ to the layer potential in an $L^{2}$ sense, and will be clarified in the error analysis at the end of the section. It should be noted that $d_{j}$ scales proportionally to a linear dimension of $\Gamma_{j}$ (for example, like $R_{j}$).

On the other hand, if
$\bx \in \mathbb{R}^{3} \setminus N_{\eta}(\Gamma_{j})$, the kernel
$G_{k}(\bx,\bX^{j}(u,v))$ in the integrand of $I^{j}_{nm}(\bx)$ is
smoother than the corresponding kernel for $\bx\in
F(\Gamma_{j})$. Thus, the above result also implies that for all
$\bx \in \mathbb{R}^{3} \setminus \Gamma_{j}$,
\begin{equation}
\sqrt{\sum_{n+m<p}\left| \tilde{I}^{j}_{nm,q}(\bx) - I^{j}_{nm} (\bx)\right|^2} \leq \varepsilon \frac{d_{j}}{\| \mtx{V} \|} \, .
\end{equation}
However, for analyzing the error in evaluating the layer potential, we wish to obtain an estimate for the contribution of a discretized patch $\Gamma_{j}$ to the layer potential $\cS_{k}[\sigma](\bx)$ denoted by $L_{j}(\bx)$,
\begin{equation}
  \label{eq:fi}
L_{j}(\bx) = \int_{u=0}^{1} \int_{v=0}^{1-u} G_{k}(\bx, \bX^{j}(u,v)) \, \sigma(u,v)
  \, J^{j}(u,v) \, du \, dv \,.
\end{equation} 
Since the density~$\sigma$ is known on patch~$\Gamma_{j}$
through its samples~$\sigma^{j}_{\ell}$ located at the
Vioreanu-Rokhlin nodes, we have that
\begin{equation}
\label{eq:sdef}
  \sigma(u,v) = \sum_{n+m < p} s^{j}_{nm} \, K_{nm}(u,v),
  \qquad \text{where} \qquad  s^{j}_{nm} = \sum_{\ell = 1}^{n_p} V_{(nm),\ell} \,
  \sigma^j_\ell ,
\end{equation}
with~$\mtx{V}$ the values-to-coefficients matrix in~\eqref{eq:v2c}.
Inserting the above expression into~\eqref{eq:fi} we have:
\begin{equation}
\label{eq:ldef}
  \begin{aligned}
L_{j}(\bx) 
  &=\sum_{n+m < p} s^{j}_{nm} 
  \int_{u=0}^{1} \int_{v=0}^{1-u} G_{k}(\bx, \bX^{j}(u,v)) \, 
    K_{nm}(u,v)  \, J^{j}(u,v) \, du \, dv \\
    &= \sum_{n+m < p} s^{j}_{nm} \, I^j_{nm}(\bx) \\
    &= \sum_{\ell = 1}^{n_p} \lp \sum_{n+m < p} V_{(nm),\ell} \,
    I^j_{nm}(\bx) \rp   \sigma^j_\ell .
  \end{aligned}
\end{equation}
Suppose next that $\tilde{L}_{j,q_{j}}(\bx)$ denotes an approximation
to $L_{j}(\bx)$, where each of the integrals $I_{nm}^{j}$ are replaced
by $\tilde{I}_{nm,q_{j}}(\bx)$. Let
$w_{j,\ell} = \sqrt{J^{j}(u_{\ell},v_{\ell}) w_{\ell}}$, and let
$\mtx{W}$ be the diagonal matrix whose entries are $w_{j,\ell}$,
$\ell=1,2\ldots n_{p}$, and let
$e_{nm} = I^{j}_{nm}(\bx) - \tilde{I}^{j}_{nm,q_{j}}(\bx)$,
$n+m<p$. Then for all
$\bx \in \mathbb{R}^{3} \setminus N_{\eta}(\Gamma_{j})$, it follows
that
\begin{equation}
\label{eq:errestfar}
\begin{aligned}
|L_{j}(\bx) - \tilde{L}_{j,q_{j}}(\bx) | &= \boldsymbol{e}^{T} \mtx{V} \mtx{W}^{-1} \mtx{W} \begin{bmatrix} \sigma^{j}_{1}  \\
\vdots \\
\sigma^{j}_{n_{p}}
\end{bmatrix} \\
& \leq \frac{\varepsilon d_{j}}{\| \mtx{V} \|} \cdot \| \mtx{V} \| \cdot \| \mtx{W}^{-1} \|  \sqrt{\sum_{\ell=1}^{n_{p}} |\sigma^{j}_{\ell}|^2 w_{j,\ell}^2}\\
& \leq \varepsilon \left( \| \sigma \|_{\mathbb{L}^{2}(\Gamma_{j})} + O(\varepsilon) \right) \, ,
\end{aligned}
\end{equation}
The last inequality follows from the fact that $\| \mtx{W}^{-1} \| = 1/d_{j}$, and that
\begin{equation}
\sum_{\ell=1}^{n_{p}} |\sigma^{j}_{\ell}|^2 w_{j,\ell}^2 = \sum_{\ell=1}^{n_{p}} |\sigma^{j}_{\ell}|^2 w_{\ell} J^{j}(u_{\ell},v_{\ell}) = \int_{\Gamma_{j}} |\sigma(\by)|^2 \, da(\by) + O(\varepsilon) \,.
\end{equation}

\begin{remark}
The same procedure as described above directly applies to the double layer potential with kernel
\mbox{$K(\bx,\by) = \bn(\by) \cdot \nabla_{\by}G_{k}(\bx,\by)$}.
For the normal derivative of the single layer potential, with kernel
\mbox{$K(\bx,\by) = \bn(\bx) \cdot \nabla_{\bx}G_{k}(\bx,\by)$}, the procedure
above can't be applied since $\bn(\bx)$ isn't well-defined at off surface
target points. However, the operator is simply the adjoint of the 
double layer, and therefore we 
use the same~$q_j$ as estimated for that case.
\end{remark}

\begin{remark}
  The $q_j$ computed in equation~\eqref{eq:qj} will of course depend
  on the kernel $G_k$, and any normalization factors. The oversampling factors were
  computed based on the Green's function $G_k(r) = e^{ikr}/(4\pi r)$.
\end{remark}

\begin{remark}
A simple calculation shows that 
\begin{equation}
L_{j,q_{j}} = \sum_{\ell=1}^{n_{q_{j}}} G_{k}(\bx, \bX^{j}(u_{\ell},v_{\ell})) J^{j}(u_{\ell},v_{\ell}) w_{\ell} \tilde{\sigma}^{j}_{\ell} \, ,
\end{equation}
where $u_{\ell},v_{\ell}$ now are the order $q_{j}$ Vioreanu-Rokhlin
nodes on $T_{0}$, $w_{\ell}$ the corresponding quadrature weights,
and $\tilde{\sigma}^{j}_{\ell}$ is the interpolated density obtained
by evaluating
\begin{equation}
\tilde{\sigma}^{j}_{\ell} = \sum_{n+m<p} s^{j}_{nm} K_{nm}(u_{\ell},v_{\ell}) \, ,
\end{equation}
where $s^{j}_{nm}$ is defined in~\cref{eq:sdef}. Since the same nodes $\bX^{j}(u_{\ell},v_{\ell})$ on $\Gamma_{j}$ can be used for all $\bx \in \mathbb{R}^{3} \setminus N_{\eta}(\Gamma_{j})$, the far part of the layer potential evaluation can be trivially coupled to fast multipole methods.
\end{remark}

\subsection{Near field quadrature}
\label{sec:nearfieldquad}

Finally, we turn our attention to the 
evaluation of~$\Snear[\sigma](\bx)$, for which the 
integrands are nearly-singular and we wish to develop a \emph{high performance}
variant of adaptive integration.
Let us consider a patch $\Gamma_{j} \in T_{\eta}(\bx)$, 
and the
integral $L_{j}(\bx)$ defined in~\cref{eq:fi}, 
which is also a near field integral for all $\bx \in N_{\eta}(\Gamma_{j})$.

It follows from~\cref{eq:ldef} that
\begin{equation}
  \begin{aligned}
  L_{j}(\bx)
    &= \sum_{\ell = 1}^{n_p} \lp \sum_{n+m < p} V_{(nm),\ell} \,
    I^j_{nm}(\bx) \rp   \sigma^j_\ell \\
    &= \sum_{\ell = 1}^{n_p} a^j_\ell(\bx) \, \sigma^j_\ell.
  \end{aligned}
\end{equation}
The numbers~$a^j_\ell(\bx)$ are the matrix entries which map the
function values~$\sigma^j_\ell$ on patch~$\Gamma_j$ to the induced
near field potential at location~$\bx$. If we approximate
each~$I^j_{nm}(\bx)$ by~$\tilde I^j_{nm}(\bx)$ to
precision~$\frac{\varepsilon d_{j}}{\| \mtx{V} \|}$, then using the same error analysis as in~\cref{sec:oversamp}, we have that
\begin{equation}
\label{eq:errestnear}
  \left|
  \sum_{\ell = 1}^{n_p} a^j_\ell(\bx) \, \sigma^j_\ell - 
  \int_{u=0}^{1} \int_{v=0}^{1-u} G_{k}(\bx, \bX^{j}(u,v)) \, \sigma(u,v)
  \, J^{j}(u,v) \, du \, dv 
  \right| \leq \varepsilon \| \sigma \|_{\mathbb{L}^{2}(\Gamma_{j})} \, .
\end{equation}

We compute~$\tilde I^j_{nm}(\bx)$ by adaptive integration on $T_{0}$. 
That is, for precision $\varepsilon$, 
we compute the integral on $T_{0}$ 
using $q$th-order Vioreanu-Rokhlin nodes, and compare it to 
the integral obtained by
\begin{enumerate}
  \item marking the midpoint of each edge of $T_0$,
  \item  subdividing $T_{0}$ into $4$ smaller right triangles, which we will call
    its {\em descendants}, and
  \item using $q$th
    order Vioreanu-Rokhlin nodes on each descendant.
\end{enumerate}
The subdivision process is repeated until, for each triangle $T$, 
its contribution to the total integral agrees with the contribution computed 
using its descendants with an error less than~$\varepsilon \cdot |T|/|T_{0}|$. 

\begin{figure}
    \centering
    \includegraphics[width=\linewidth]{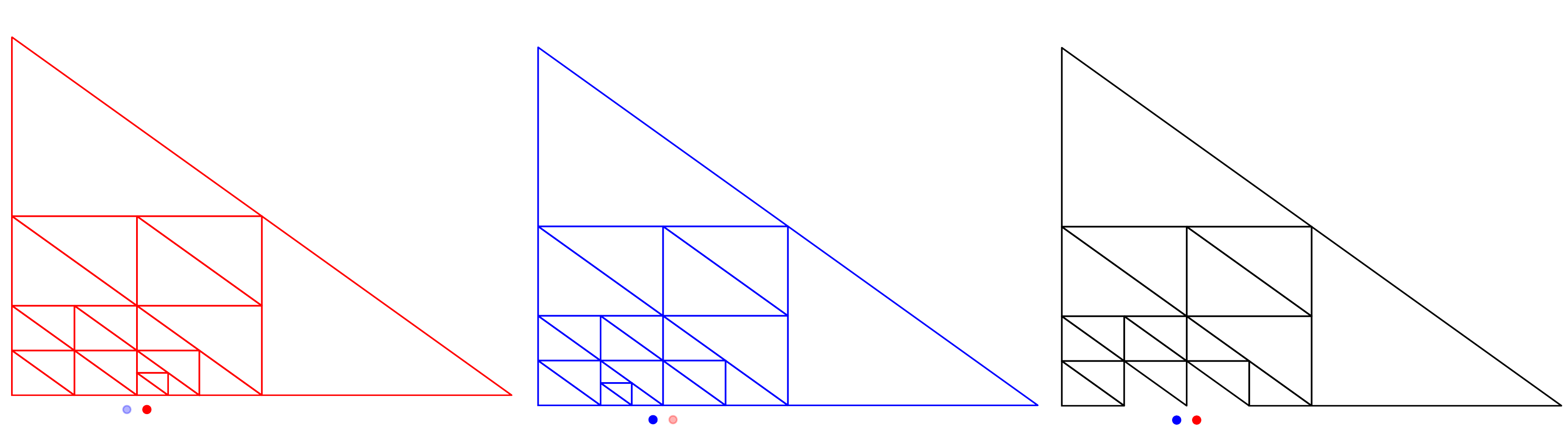}
    \caption{Adaptive integration grids used for the red target (left), blue target (center). The black grid (right) is the common set of triangles in both of the grids for which the function values of $\bX^{j},K_{nm}$, and $J^{j}$ are reutilized.}
    \label{fig:multipletargpot}
\end{figure}

Done naively, this adaptive integration process dominates the cost
of quadrature generation because of the large number of targets 
in $N_{\eta}(\Gamma_{j})$. Note, however, that
as we vary~$n,m$, for a fixed target~$\bx$, the integrand 
of $I^j_{nm}(\bx)$ includes the same kernel values~$G(\bx,\bX^{j}(u,v))$.
Moreover, the adaptive grids generated for different targets 
have significant commonality. 
Thus, we can  reuse the function values
of~$\bX^{j}(u,v)$,~$K_{nm}(u,v)$
and~$J^{j}(u,v)$ if they have already been computed on any descendant
triangle (see Fig.~\ref{fig:multipletargpot}).
The resulting scheme incurs very little increase in storage 
requirements; this significantly improves the overall performance. 

\begin{remark}
Adaptive integration often results in much greater accuracy 
than requested.
With this in mind, we set $\varepsilon$ in the termination
criterion to be somewhat larger 
than the precision requested. Our choice is based on
extensive numerical experimentation and 
the full set of parameters used in our implementation is
 available at~\url{https://gitlab.com/fastalgorithms/fmm3dbie}.
\end{remark}

\begin{remark}
To further improve the performance of computing~$I^{j}_{nm}(\bx)$, we 
make use of {\em two} parameters:~$\eta$ and~$\eta_1 < \eta$.
We only use adaptive integration for the nearest targets, inside
$N_{\eta_1}(\Gamma_j)$. For targets
$\bx \in N_{\eta}(\Gamma_{j}) \setminus N_{\eta_{1}}(\Gamma_{j}) $, we
use a single oversampled quadrature without any adaptivity.
This is slightly more expensive in terms of function evaluations,
but eliminates the branching queries of adaptive quadrature and allows the 
use of highly optimized linear algebra libraries. 
From extensive numerical experiments, we have found that
$\eta_{1} = 1.25$ provides a significant speedup.
\end{remark}

\subsection{Error Analysis}
There are two sources of error in the computation of
$\cS_{k}[\sigma]$: the discretization error due to discretizing the
charts $\bX^{j}$ using order $p$ Vioreanu-Rokhlin nodes, and the
quadrature error due to using different quadrature rules for
evaluating integrals over the discretized patches. As shown
in~\cite{Atkinson95}, the discretization error can be bounded by
\begin{multline}
\label{eq:errdisc}
\left| \cS_{k}[\sigma](\bx) - \sum_{j=1}^{\Npat}  \int_{u=0}^{1} \int_{v=0}^{1-u} G_{k}(\bx, \bX^{j}(u,v)) \sigma(u,v) J^{j}(u,v) \, du \, dv \right| \\
= \left| \cS_{k}[\sigma](\bx) - \sum_{j=1}^{\Npat} L_{j}(\bx) \right| \leq C h^{p},
\end{multline}
where $h = \max_{j} R_{j}$, and some domain dependent constant
$C$. For a
given $\bx \in \Gamma_{m}$, our method approximates the layer potential as 
\begin{multline}
I(\bx) = \sum_{\ell=1}^{N(\bx)} G_{k}(\bx, \bX^{m}(u_{\bx,\ell},v_{\bx,\ell})) \, \sigma(u_{\bx,\ell},v_{\bx,\ell})
  \, J^{m}(u_{\bx,\ell},v_{\bx,\ell}) w_{\bx,\ell} \  + \\ 
 \sum_{\Gamma_{j} \in N_{\eta}(\bx)} \sum_{\ell=1}^{n_{p}}
   a_{\ell}^{j}(\bx) \sigma^{j}_{\ell}  \ + 
 \sum_{\Gamma_{j} \in F_{\eta}(\bx)} \sum_{\ell=1}^{n_{q_{j}}}
   G_{k}(\bx, \bX^{j}(u_{\ell,q_{j}},v_{\ell,q_{j}}))
   J^{j}(u_{\ell,q_{j}},v_{\ell,q_{j}}) w_{\ell,q_{j}}
   \tilde{\sigma}^{j}_{\ell}  \, ,
\end{multline}
which are approximations to $\Sself$, $\Snear$, and $\Sfar$
respectively. Here $(u_{\bx,\ell},v_{\bx,\ell}), w_{\bx,\ell}$ are the
auxiliary nodes on the self patch, $a^{j}_{\ell}(\bx)$ are the near
quadrature corrections computed via adpative integration,
$\tilde{\sigma}^{j}_{\ell}$ is the oversampled density, and
$(u_{\ell,q},v_{\ell,q}), w_{\ell,q}$, $\ell=1,2,\ldots n_{q}$ are the
order $q$ Vioreanu-Rokhlin nodes on $T_{0}$.  Using the estimates
in~\cref{eq:errestself,eq:errestnear,eq:errestfar}, combined
with~\cref{eq:errdisc}, we conclude that
\begin{equation}
\left| \cS_{k}[\sigma](\bx) - I(\bx) \right| \leq C h^{p} + \varepsilon \| \sigma \|_{\mathbb{L}^{2}(\Gamma)} \,.
\end{equation}
\begin{remark}
  As shown in~\cite{Atkinson95}, the discretization error for
  evaluating $\cD_{k}[\sigma]$ is $O(h^{p-1})$. The quadrature error
  analysis remains the same and we can evaluate $\cD_{k}[\sigma](\bx)$
  with accuracy
  $O(h^{p-1}) + \varepsilon \| \sigma \|_{\mathbb{L}^{2}(\Gamma)}$.
\end{remark}

\section{Coupling quadratures to FMMs}
\label{sec:fmmcoupling}

For a complete description of three-dimensional
FMMs applied to sums of the form \eqref{fmmpt}, we refer the reader to the
original papers \cite{fmm2,wideband3d,greengard-1997,greengard-huang}.
In order to understand the modifications needed for evaluating layer 
potentials, however, we will need to make reference to the adaptive oct-tree
data structures on which the FMM is built. We briefly summarize that
construction here.

\subsection{Level-restricted, adaptive oct-trees}

Suppose for the moment that we are given a collection of $N$ points,
contained in a cube $C$.  We will superimpose on $C$ a hierarchy of 
refinements as follows:
the root of the tree is $C$ itself and defined as {\em level 0}.
Level $l+1$ is obtained from level
$l$ recursively by subdividing each cube at level $l$ into eight equal parts,
so long as the number of points in that cube at level $l$ is greater than
some specified parameter $s$.
 The eight cubes created in the above step are referred to as its 
children.  Conversely, the box which was divided is referred to as their
parent.
When the refinement has terminated, $C$ is covered
by disjoint childless boxes at various levels of the hierarchy (depending
on the local density of the given points). These childless boxes are referred
to as leaf nodes.  For any box $D$ in the hierarchy, 
other boxes at the same level that touch $D$ are called its {\em
  colleagues}.
For simplicity, 
we assume that the oct-tree satisfies a standard restriction - namely, 
that two leaf nodes which share a boundary point must be no more than one 
refinement level apart. In creating the adaptive data structure as described
above, it is very likely that the level-restriction criterion is
not met. Fortunately, assuming that the 
tree constructed to this point has $O(N)$ leaf nodes and that its depth is 
of the order $O(\log N)$, it is straightforward to enforce the
level-restriction in a second step requiring~$O(N \log N)$ effort 
 with only a modest amount of additional refinement~\cite{treebook}.

\subsection{Precomputation}

To reiterate, on input, we assume we are given  
a surface~$\Gamma$ consisting of (curvilinear) triangles~$\Gamma_j$,
\begin{equation}
  \Gamma = \cup_{j=1}^{\Npat} \, \Gamma_{j},
\end{equation}
each given to the desired order of accuracy~$p$. Each~$\Gamma_j$ is
then discretized using~$n_p$ points which are the images under the
map~$\bX^j: T_0 \to \bbR^3$ of the Vioreanu-Rokhlin nodes on the
standard simplex~$T_0$.  We will refer to these as the {\em
  discretization nodes}, on which we assume that samples of the
density~$\sigma$ are known.  The total number of such points is $N =
\Npat \times n_p$.  As above, we let~$\bc_j$ denote the centroid of
the $j$th patch and $R_j$ the radius of the smallest sphere centered
at $\bc_j$ that contains~$\Gamma_j$.  We assume there are $N_T$
targets, which could be either the discretization nodes themselves, a
collection of off-surface points, or both.

In coupling the FMM to local quadratures, we need to determine, for
each surface patch, which targets are in its near field and what order
Vioreanu-Rokhlin quadratures are needed for the far field
computation.  Both are controlled by the parameter $\eta$, as
discussed in \Cref{sec:etachoice}. The default value for $\eta$ is
$2.75$, $2$, or $1.25$ depending on whether the desired order of
accuracy is $p \leq 4$, $4 < p \leq 8$ or $p> 8$,
respectively.  The
first step is to build an adaptive oct-tree based on the patch
centroids $\{ \bc_j \}$ and the target locations, with one minor
modification.  That modification is to prevent triangle centroids
associated with large triangles from propagating to fine levels during
the tree construction.  For this, suppose $\bc_j$ is in some box,
denoted $D(\bc_j)$ at level $l$, and let $d$ denote the linear
dimension of $D(\bc_j)$. If $2 \eta \, R_j > d$,
then we leave the centroid associated with $D(\bc_j)$, while allowing
smaller triangles and/or targets to be associated with the
children. We will say that $\Gamma_j$ is {\em tethered} at level $l$.

The near field for each patch is now easy to determine.  For each
triangle $\Gamma_j$, let $D(\bc_j)$ denote the box to which the
triangle centroid is associated - either a leaf node or the box at a
coarser level $l$ if is tethered there.  Clearly, if $D(\bc_{j})$ is
not a leaf box, then the near field region $B_{\eta R_j}(\bc_j)$ is
contained within $D(\bc_j)$ and its colleagues. If $D(\bc_{j})$
  is a leaf box, then the near field region $B_{\eta R_j}(\bc_j)$ is
  contained within the union of $D(\bc_{j})$, its colleagues, and leaf
  boxes which are larger in size than $D(\bc_{j})$ and which share a
  boundary with $D(\bc_{j})$.  Scanning those colleagues, all targets
$\bx$ that do not lie on $\Gamma_j$ itself and satisfy the criterion
\[ |\bx-\bc_{j}| < \eta R_{j}     \]
are assigned to the {\em near field list} for $\Gamma_j$.
One can then compute the near field quadratures using the method of 
\Cref{sec:nearfieldquad} for each point in the target
list. 
This requires storing a matrix of dimension 
$N_{near}(j) \times n_{p}$, where
$N_{near}(j)$ is the size of the target list.
We will denote this matrix by ${\cal N}_j$.

Assuming one wishes to evaluate the layer potential on surface,
we also need to compute 
the self interactions for each triangle
using the generalized Gaussian quadrature scheme of 
\cite{bremer_2012c,bremer_2013}, as described in \Cref{sec:local}.
This requires storing an $n_p \times n_p$ matrix
for each patch, which we will denote by 
${\cal S}_j$. 

Once the near field work has been carried out, 
the far field quadrature
order $q_j$ is determined, as described in \Cref{sec:oversamp}.
One can then interpolate from the $n_p$ discretization nodes 
on $\Gamma_j$ to the $n_{q_j}$ quadrature nodes on $\Gamma_j$ 
using the Koornwinder basis for interpolation. We will denote by 
$N_{over}$ the total number of oversampled points used:
$N_{over} = \sum_{i =
    1}^{\Npat} n_{q_i}$. 

\subsection{Fast evaluation of layer potentials}

The simplest FMM-based scheme for evaluating a layer potential
is to call the point-based FMM in the form \eqref{fmmpt},
with $N_{over}$ sources and $N_T$ targets. 
For every target, if it is in the near field of 
patch $\Gamma_j$, one subtracts the contribution made in the naive,
point-based FMM calculation from the $n_q$ oversampled points on that patch.
The potential at the target can then be incremented by the appropriate,
near field quadrature-corrected
interactions, using the stored matrix ${\cal N}_j$. 
If the target is on surface (one of the 
discretization nodes on $\Gamma_j$ itself),
the correct self interaction is obtained from the 
precomputed matrix ${\cal S}_j$.

We refer to the algorithm above as the
\emph{subtract-and-add} method.
It has the drawback that it could suffer from
catastrophic cancellation for dense discretizations with 
highly adaptive oct-trees since the near field point contributions
within the naive FMM call are spurious and could be much greater in magnitude
than the correct contributions.
(In practice, we have not detected any such loss of accuracy, at least  
for single or double layer potentials.)

For readers familiar with the FMM, it is clear that one could avoid the
need to compute and then subtract spurious contributions, by disabling
the direct (near neighbor) interaction step in the FMM.
When looping over the leaf nodes, for each source-target pair, one can first 
determine whether the source is on a patch for which the target
is in the far field. If it is, carry out the direct interaction. If it is 
in the near field, omit the direct interaction.
When the FMM step is completed, the subsequent processing
takes place as before, but there is no need to subtract any spurious 
contributions.

It is, perhaps, surprising that the \emph{subtract-and-add} method is 
faster in our current implementation, 
even though more flops are executed. This is
largely because of the logical overhead and the bottlenecks
introduced in loop unrolling and other compiler-level code optimizations.

\begin{remark}
The adaptive oct-tree used in the point FMM is different from the one
used for determining the near field of the patches. The latter is
constructed based on centroid and target locations, while the point
FMM oct-tree is constructed based on \emph{oversampled} source and
target locations. Thus, different termination criteria can be chosen
for the construction of these oct-trees in order to optimize the
performance of the separate tasks.

Since the additional processing required for evaluating layer
potentials is decoupled from the algorithm used for accelerating the
far field interactions, one could use any fast hierarchical algorithm
like the FMM, an FFT-based scheme like fast Ewald summation, or a
multigrid-type PDE solver.
\end{remark}

\section{Numerical examples}
\label{sec:examples}
\begin{figure}
    \centering
    \includegraphics[width=0.6\linewidth]{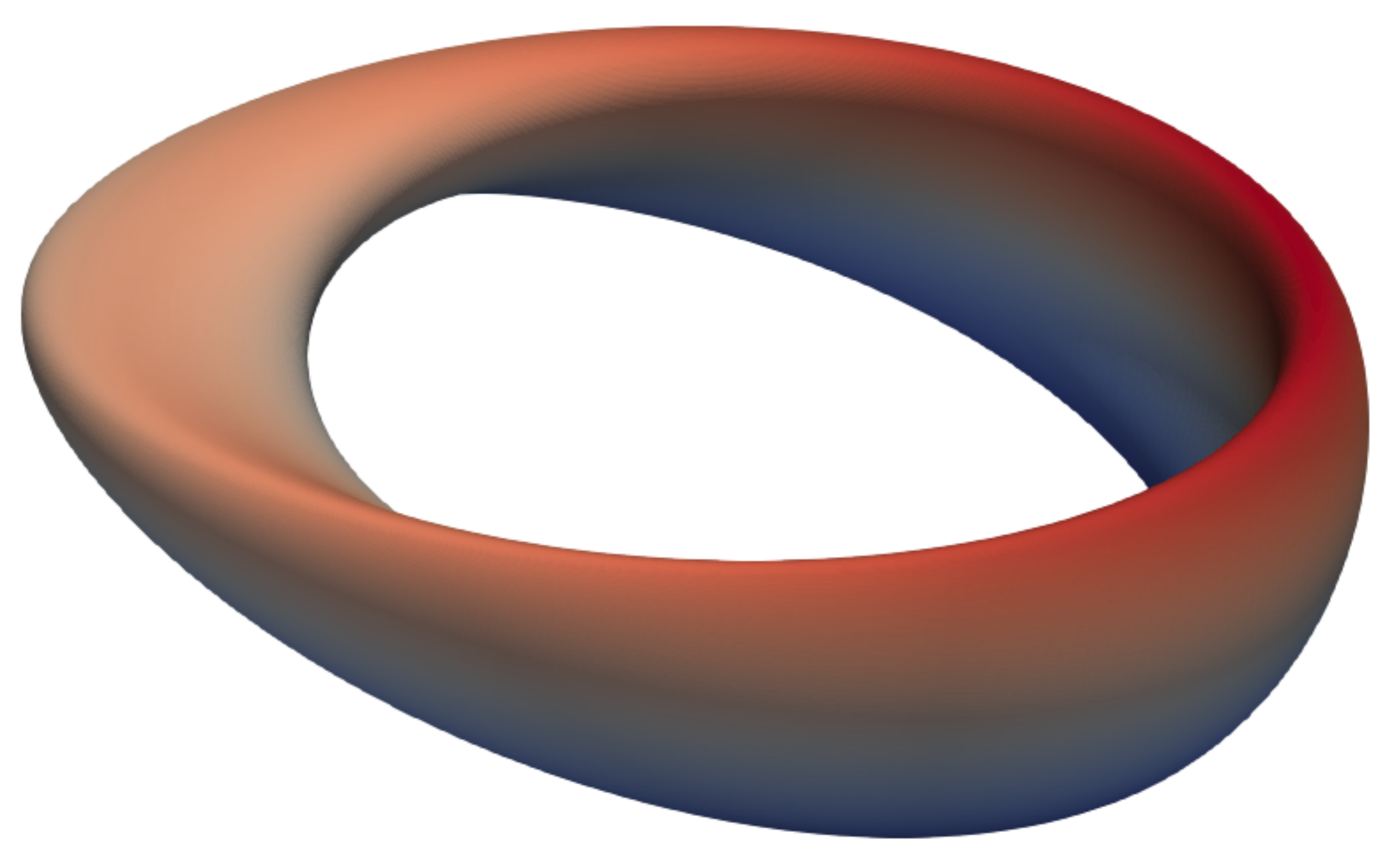}
    \caption{The boundary of a stellarator like geometry. The surface is colored using its z-coordinate}
    \label{fig:stell-disc}
\end{figure}

In this section, we illustrate the performance of our approach.
For Examples \ref{subsec-cpu}, \ref{subsec-aspect-ratio}, and
\ref{subsec-num-order}, we consider a twisted torus as the geometry
(typical of stellarator design in plasma physics applications).
The boundary $\Gamma$ is parameterized by $\bX: [0,2\pi]^2 \to \Gamma$
with
\begin{equation}
    \bX(u,v) = \sum_{i=-1}^{2} \sum_{j=-1}^{1} \delta_{i,j} \begin{bmatrix} \cos{v} \cos{((1-i)\, u+j\, v)} \\
    \sin{v} \cos{((1-i) \, u+j \, v)} \\
    \sin{((1-i)\, u+j \, v)} 
    \end{bmatrix} \, ,
\end{equation}
where the non-zero coefficients are $\delta_{-1,-1}=0.17$, 
$\delta_{-1,0} = 0.11$, $\delta_{0,0}=1$, $\delta_{1,0}=4.5$, $\delta_{2,0}=-0.25$, $\delta_{0,1} = 0.07$, and $\delta_{2,1}= -0.45$.
(See Fig. \ref{fig:stell-disc}.) 

The code was implemented in Fortran and 
compiled using the GNU Fortran 10.2.0 compiler.
We use the point-based FMMs from the FMM3D package (\url{https://github.com/flatironinstitute/FMM3D}). All CPU timings in these examples were obtained on a laptop
using a single core of an Intel i5 2.3~GHz processor. 

We use the following metrics to demonstrate the performance of our approach.
As above, for discretization order $p$, we let 
$n_p = p(p+1)/2$, and we let $q_{j}$ denote the far-field quadrature order 
for $\Gamma_{j}$. 
The user-specified precision is denoted by $\varepsilon$.
Recall that the total
number of discretization points 
on the boundary is denoted by $N = \Npat \cdot n_p$ and that 
the total number of oversampled nodes
is denoted by $\Nover = \sum_{j=1}^{\Npat} n_{q_j}$.
We define the oversampling parameter by $\alpha= \Nover/N$.
The memory requirements per discretization node for storing all
interactions in $\Snear$ are given by
\begin{equation}
m = \frac{n_p \left( \sum_{j=1}^{\Npat}  N_{near}(j) + n_p \right)}{N}.
\end{equation}
This accounts for both off-surface targets and on surface evaluation.

Let $\tinit$ denote the time required to precompute all
near field quadrature corrections and let $\tlp$ denote the time for 
evaluating the layer potential given the precomputed near field quadratures 
Then, the quantities $\sinit= N/\tinit$, and $\slp=N/\tlp$, 
are the speeds of the corresponding steps,
measured in points processed per second.

One feature of the surface triangulation that has some influence on
speed is the aspect ratio of the patches.
Letting $\sigma_{1},\sigma_{2}$ be the eigenvalues of the first 
fundamental form of $\Gamma_{j}$, we define its aspect ratio by
\begin{equation}
a_{j} = \sqrt{\frac{\int_{\Gamma_{j}} 
\left(\frac{\sigma_{1}}{\sigma_{2}}\right)^2 
da(\by)}{\int_{\Gamma_{j}} da(\by)}} \, .
\end{equation}
We let $\amax = \max_{j} a_{j}$ and $\aavg = \sum_{j} a_{j}/\Npat$,
the maximum aspect ratio and the average aspect ratio over all
triangles, respectively.

\subsection{Memory and oversampling requirements \label{subsec-cpu}}

To illustrate the performance of our method as a function of the order of
accuracy $p$ and the requested precision $\varepsilon$, we 
consider the evaluation of the single layer potential
$\cS[\sigma]$ with frequency $k=1$ on the stellarator geometry 
discretized with $\Npat=2400$ (the diameter of the stellarator with $k=1$ is approximately 1.7 wavelengths).
In~\Cref{tab:numerical-perf-res}, we tabulate the memory requirements
per point $m$, the oversampling $\alpha$, and the speeds 
$\sinit$ and $\slp$, as we vary $p$ and $\varepsilon$. 
The scheme behaves as expected: for fixed $p$, as $\varepsilon \rightarrow 0$, 
$\alpha$ increases while $\sinit$ and $\slp$ decrease. The oversampling 
parameter $\alpha$ depends on both $p$ and $\varepsilon$, as discussed in 
\Cref{sec:oversamp}.

\afterpage{
\begin{table}
\begin{minipage}{0.5\linewidth}
\begin{center}
\[\begin{array}{|c|c|c|c|c|}
\hline\tikz{
\node[below left, inner sep=1pt] (v1) {$p$};
\node[above right, inner sep=1pt] (v2) {$\varepsilon$};
\draw(v1.north west|-v2.north west) -- (v1.south east-|v2.south east);}
 & 5 \cdot 10^{-3} & 5 \cdot 10^{-4} & 5 \cdot 10^{-7} & 5 \cdot 10^{-10}\\ \hline 
2 & 1 & 1 & 3.52 & 9.35\\ \hline 
3 & 0.5 & 1 & 2.5 & 5.19\\ \hline 
4 & 0.6 & 0.824 & 2.8 & 4.5\\ \hline 
6 & 0.476 & 0.714 & 1.84 & 2.94\\ \hline 
8 & 0.777 & 0.923 & 2.24 & 3.9\\ \hline 
\end{array}\]
\subcaption{Oversampling parameter $\alpha$}
\end{center}
\end{minipage}
\begin{minipage}{0.5\linewidth}
\begin{center}
\[\begin{array}{|c|c|c|c|c|}
\hline\tikz{
\node[below left, inner sep=1pt] (v1) {$p$};
\node[above right, inner sep=1pt] (v2) {$\varepsilon$};
\draw(v1.north west|-v2.north west) -- (v1.south east-|v2.south east);}
 & 5 \cdot 10^{-3} & 5 \cdot 10^{-4} & 5 \cdot 10^{-7} & 5 \cdot 10^{-10}\\ \hline 
2 & 127 & 127 & 127 & 127\\ \hline 
3 & 258 & 258 & 258 & 258\\ \hline 
4 & 209 & 209 & 209 & 209\\ \hline 
6 & 440 & 440 & 440 & 440\\ \hline 
8 & 286 & 286 & 286 & 286\\ \hline 
\end{array}\]
\subcaption{Memory requirements $m$ per point }
\end{center}
\end{minipage}
\newline
\begin{minipage}{0.5\linewidth}
\begin{center}
\[\begin{array}{|c|c|c|c|c|}
\hline\tikz{
\node[below left, inner sep=1pt] (v1) {$p$};
\node[above right, inner sep=1pt] (v2) {$\varepsilon$};
\draw(v1.north west|-v2.north west) -- (v1.south east-|v2.south east);}
 & 5 \cdot 10^{-3} & 5 \cdot 10^{-4} & 5 \cdot 10^{-7} & 5 \cdot 10^{-10}\\ \hline 
2 & 8460 & 6480 & 2970 & 1350\\ \hline 
3 & 8130 & 5880 & 2460 & 1150\\ \hline 
4 & 8230 & 6440 & 2820 & 1360\\ \hline 
6 & 2890 & 2380 & 1250 & 734\\ \hline 
8 & 2080 & 1750 & 945 & 581\\ \hline 
\end{array}\]
\subcaption{Precomputation speed $\sinit$}
\end{center}
\end{minipage}
\begin{minipage}{0.5\linewidth}
\begin{center}
\[\begin{array}{|c|c|c|c|c|}
\hline\tikz{
\node[below left, inner sep=1pt] (v1) {$p$};
\node[above right, inner sep=1pt] (v2) {$\varepsilon$};
\draw(v1.north west|-v2.north west) -- (v1.south east-|v2.south east);}
 & 5 \cdot 10^{-3} & 5 \cdot 10^{-4} & 5 \cdot 10^{-7} & 5 \cdot 10^{-10}\\ \hline 
2 & 30100 & 24500 & 6480 & 1450\\ \hline 
3 & 35800 & 18600 & 6870 & 2170\\ \hline 
4 & 27500 & 20200 & 5320 & 2350\\ \hline 
6 & 24300 & 19900 & 6290 & 3390\\ \hline 
8 & 20500 & 16000 & 5960 & 2140\\ \hline 
\end{array}\]
\subcaption{Layer potential speed $\slp$}
\end{center}
\end{minipage}
\caption{
Memory requirements, oversampling requirements, precomputation speed and
layer potential evaluation speed
as a function of order of accuracy $p$ and precision $\varepsilon$}
\label{tab:numerical-perf-res}
\end{table}

\begin{table}
\begin{center}
\[\begin{array}{|c|c|c|c|c|c|}
\hline
\aavg & 1.61 & 2.33 & 4.66 & 9.32 & 14\\ \hline
\Npat & 2400 & 2304 & 2450 & 2304 & 2400\\ \hline 
\alpha & 2.8 & 2.8 & 2.8 & 2.8 & 2.8\\ \hline 
m & 209 & 291 & 585 & 1220 & 1860\\ \hline 
\sinit & 2890 & 2160 & 1100 & 519 & 341\\ \hline 
\squad & 5390 & 4940 & 3920 & 2890 & 2270\\
\hline 
\end{array}\]
\end{center}
\caption{Performance as a function of average aspect ratio $\aavg$}
\label{tab:numerical-asp-res}
\end{table}

\begin{figure}
  \centering
  \begin{subfigure}[t]{.45\linewidth}
    \centering
    \includegraphics[width=.95\linewidth]{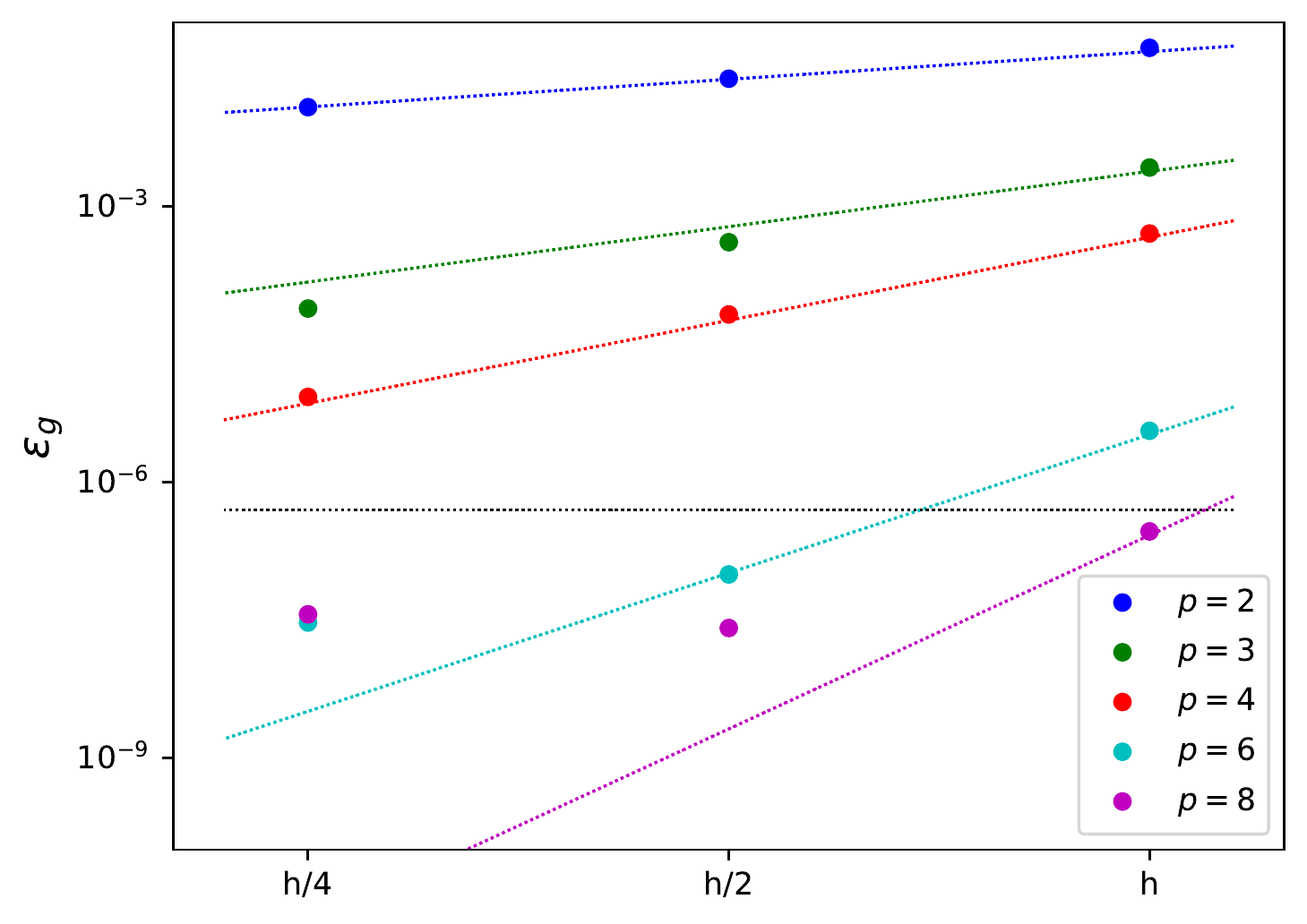}
    \caption{Relative $L^2$ error in Green's identity, denoted by
  $\varepsilon_{g}$.}
  \end{subfigure}
  \quad
  \begin{subfigure}[t]{.45\linewidth}
    \centering
    \includegraphics[width=.95\linewidth]{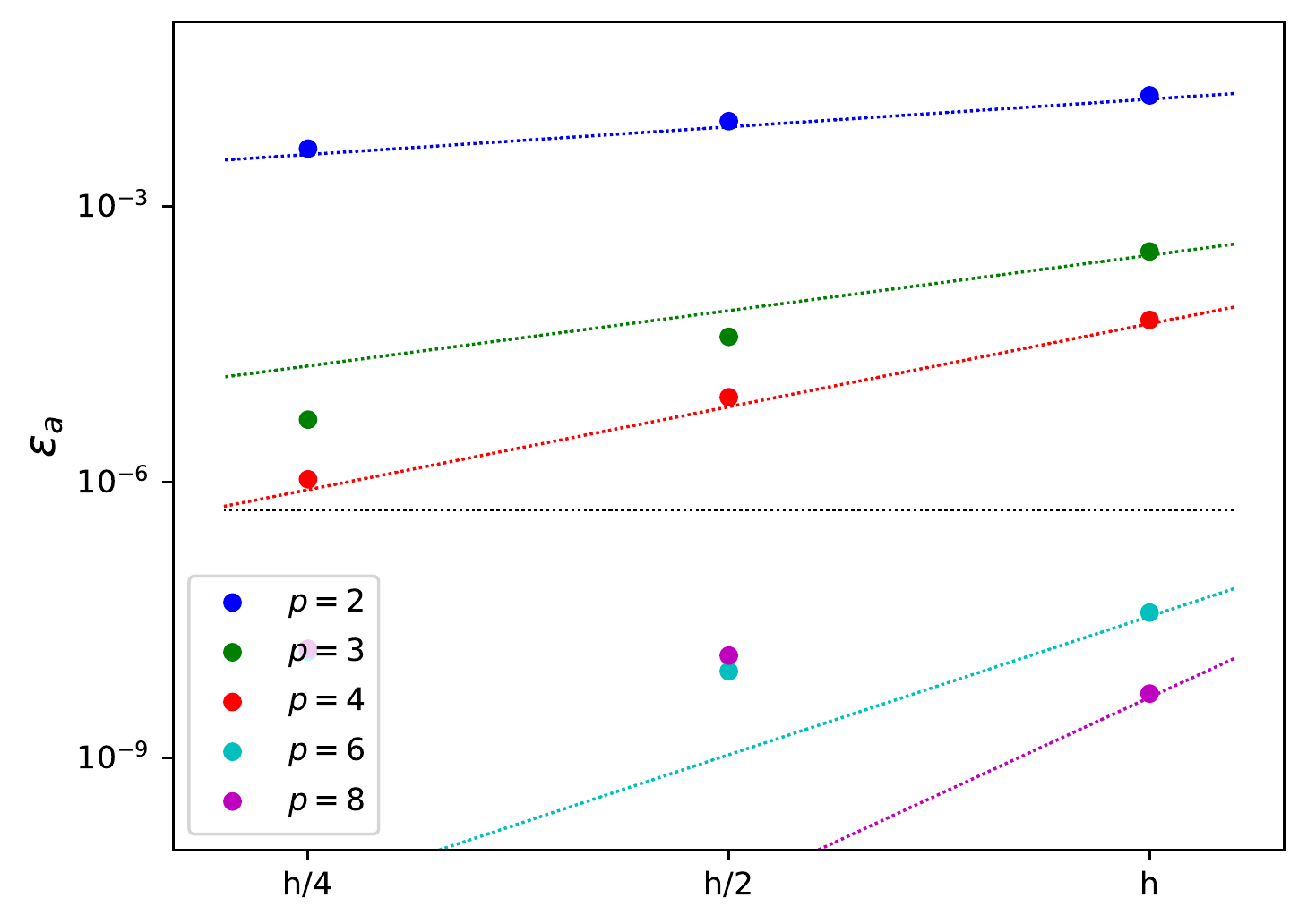}
    \caption{Relative $L^{\infty}$ error in solution to
  integral equation, denoted by  $\varepsilon_{a}$.}
  \end{subfigure}
\caption{Relative errors of layer potential evaluations.  In both
  figures, the dashed colored lines are reference curves for $h^{p-1}$
  with the corresponding $p$. The dashed black line is a reference
  line for the specified tolerance~$\varepsilon$. }
\label{fig:ooc}
\end{figure}
\clearpage
}

\subsection{Effect of aspect ratio \label{subsec-aspect-ratio}}

To investigate the effect of triangle quality on the performance of our
method, we vary the average aspect ratio of the discretization. The task 
at hand is again to compute $\cS[\sigma]$ with $k=1$ on the stellarator, 
while varying the triangulation without 
a significant change in the total number of patches.
In~\Cref{tab:numerical-asp-res}, we tabulate $\aavg$, $\Npat$, $\alpha$, 
$m$, $\sinit$, and $\slp$ for $p=4$ and $\varepsilon = 5 \cdot 10^{-7}$. 
We note that (except for the oversampling parameter), the performance of 
the approach deteriorates as the average aspect ratio of the 
discretization is increased, especially in the precomputation phase.

\subsection{Order of convergence \label{subsec-num-order}}

To demonstrate the accuracy of our approach, we consider two tests.
First, we verify Green's identity along the surface:
\begin{equation}
    \frac{u}{2} = \cS_{k} \left[\frac{\partial u}{\partial n}\right] - \cD_{k}[u] \, , 
\end{equation}
where $u$ is the solution to the Helmholtz equation in the interior of the 
domain $\Omega$ generated by a point source located in the exterior. 
The second test is to use the combined field representation
\eqref{eq:urep} to solve the Dirichlet problem for an unknown density~$\sigma$.
With the right-hand side in the corresponding integral equation 
\eqref{eq:inteq} taken to be a known solution $u$, $\sigma$ satisfies
\begin{equation}
    \frac{\sigma}{2} + \cD_{k}[\sigma] - ik \cS_{k}[\sigma] = u,
    \qquad \text{on } \Gamma.
\end{equation}
We can then check that $u$ is correctly reproduced at any point in the
exterior.  For both of these tests, we discretize the stellarator
using $\Npat = 600$, $2400$, and $9600$, and compute the layer
potentials with a tolerance of $\varepsilon = 5 \times
10^{-7}$ for both the quadratures and the FMM.

In~\Cref{fig:ooc}, we plot the relative $L^{2}$ error in Green's 
identity $\varepsilon_{g}$ (left), and the relative $L^{\infty}$ error 
at a point in the interior $\varepsilon_{a}$ (right) as we vary the order 
of discretization. Note that in both tests, the errors decrease at the
rate $h^{p-1}$ until the tolerance $\varepsilon$ is reached.
This is consistent with the analysis in \cite{Atkinson95}.

\subsection{Large-scale examples}
We demonstrate the performance of our solver on several large-scale
problems.  We first solve for the electrostatic field induced by an
interdigitated capacitor, followed by Dirichlet and Neumann boundary 
value problems governed by the Helmholtz equation in the exterior of
an aircraft. Our last example involves scattering in a medium with 
multiple sound speeds, modeled after a Fresnel lens. The results in this
section were obtained using an Intel Xeon Gold 6128 Desktop with 24
cores.

\subsubsection*{Interdigitated capacitor}
\label{subsec:cap}

\begin{figure}[t]
\centering
  \includegraphics[width=0.9\linewidth]{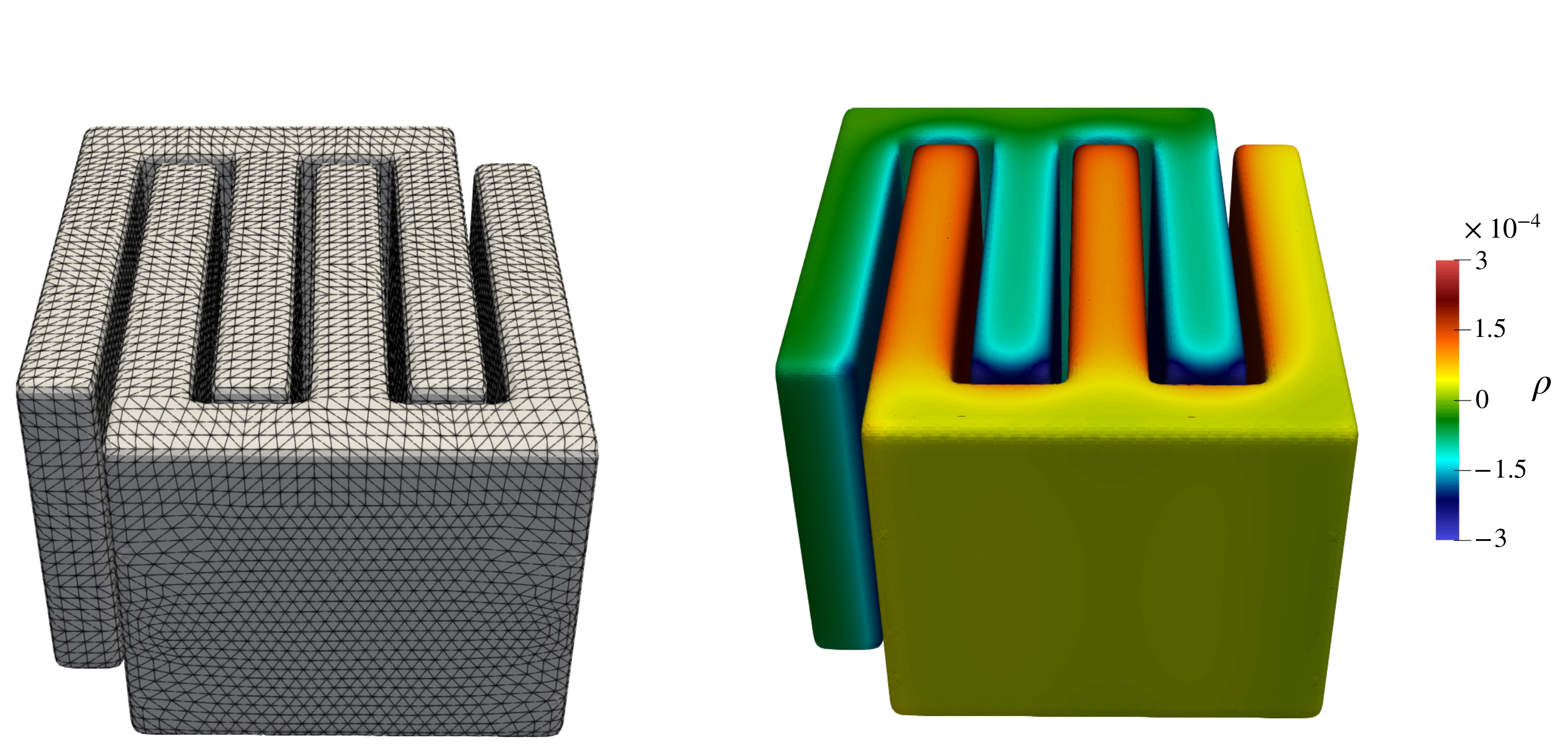}
  \caption{(left) A 4th-order computational mesh of the boundary, and (right) the solution $\rho$ to~\cref{eq:cap-inteq}}
  \label{fig:cap}
\end{figure}

A challenging problem in electrostatics is the calculation of the 
capacitance of a configuration of two perfect compact conductors with 
complicated contours, which may also be close to touching. 
(See~\cref{fig:cap}.)
The capacitance is defined as the ratio $C = Q/V$, where 
$V$ is the potential difference between the conductors, $Q$ is the 
total charge held on one conductor and $-Q$ is the total charge
held on the other.
In simulations, $C$ can be computed in two ways.
First, one can solve the Dirichlet problem for 
the electrostatic potential
$u$, with $u= 0$ on one conductor and $u = 1$ on the other.
From the computed solution, the total charge can be obtained via
the integral~\cite{jackson}
\begin{equation}
  Q = \int_\Gamma \frac{\partial u}{\partial n}(\bx')  \, da(\bx').
\end{equation}
The capacitance is then $C = Q/1 = Q$.

A second (equivalent) approach, which we will take here, is to place a
net charge~$q_{1} = -1$ on one conductor~$\Omega_1$ with boundary
$\Gamma_{1}$ and a net charge of~$q_{2} = 1$ on the other
conductor~$\Omega_2$ with boundary $\Gamma_{2}$.  One can then
determine the corresponding potential difference by solving the
following boundary value problem for the potential~$u$ in the
domain~$E$ exterior to $\Omega_1$ and $\Omega_2$, i.e. the
domain~$E = \bbR^3 \setminus (\Omega_1 \cup \Omega_2)$:
\begin{equation}
  \begin{aligned}
    \Delta u &= 0,  &\quad &\bx \in E , \\
    u &= V_{i}, & &\bx \text{ on } \Gamma_{i}, , \\
-\int_{\Gamma_{j}} \frac{\partial u}{\partial n} \, da &= q_{i}, & & \\
u &\to 0 & &\text{ as } |\bx| \to \infty.
\end{aligned}
\end{equation}
Here, the constants $V_{1}$, $V_{2}$ are unknowns as well as the
potential $u$. The \emph{elastance} of the system is then given by $P
= (V_{2}-V_{1})/(q_{2}-q_{1}) = (V_{2}-V_{1})/2$. It is the inverse of
the corresponding capacitance $C = 1/P = 2/(V_{2}-V_{1})$.  This
formulation (which can involve more than two conductors) is generally
referred to as the \emph{elastance} problem (see~\cite{elastance-ref}
and the references therein). PDEs of this type where Dirichlet data is
specified up to an unknown constant are sometimes called modified
Dirichlet problems~\cite{mikhlin-book}.

We represent the solution $u$ using the combined field representation, 
\begin{equation}
u = S_{0}[\rho] + D_{0}[\rho] ,
\end{equation}
where~$\rho$ is an unknown density. Imposing the boundary conditions 
on~$\Gamma_{i}$, $\rho$ must satisfy the integral equation
\begin{equation}
  \begin{aligned}
    \rho/2 + D_{0}[\rho] + S_{0}[\rho] &= V_{i},   &\quad
        &\bx \text{ on } \Gamma_{i} , \\
    \int_{\Gamma_{i}} \rho &= -q_{i}, & &i=1,2 .
\end{aligned}
\label{eq:cap-inteq}
\end{equation}
We discretize the surface~$\Gamma$ with~$\Npat=29,888$ and~$p=4$, and then
solve the resulting linear system of size $N=298,882$ using GMRES. We 
set the quadrature tolerance $\varepsilon=5\times 10^{-7}$. For this 
setup, $\tinit=6.7$s, $\alpha=2.98$, $m=154.2$. GMRES converged to 
a relative residual of $5 \times 10^{-7}$ in 25 iterations, and the 
solution was obtained in $128.4$s.  The reference capacitance for the 
system was computed by refining each patch until it had converged to
5 significant digits, given by $2237.1$. The relative error in the 
computed capacitance was $2.2 \times 10^{-4}$. In~\cref{fig:cap}, 
we plot the computational mesh and the solution $\rho$ on the surface of the conductors.

\subsubsection*{Scattering from an airplane (sound-soft)}
\label{subsec:plane}

\begin{figure}[b]
\begin{center}
    \includegraphics[width=0.75\linewidth]{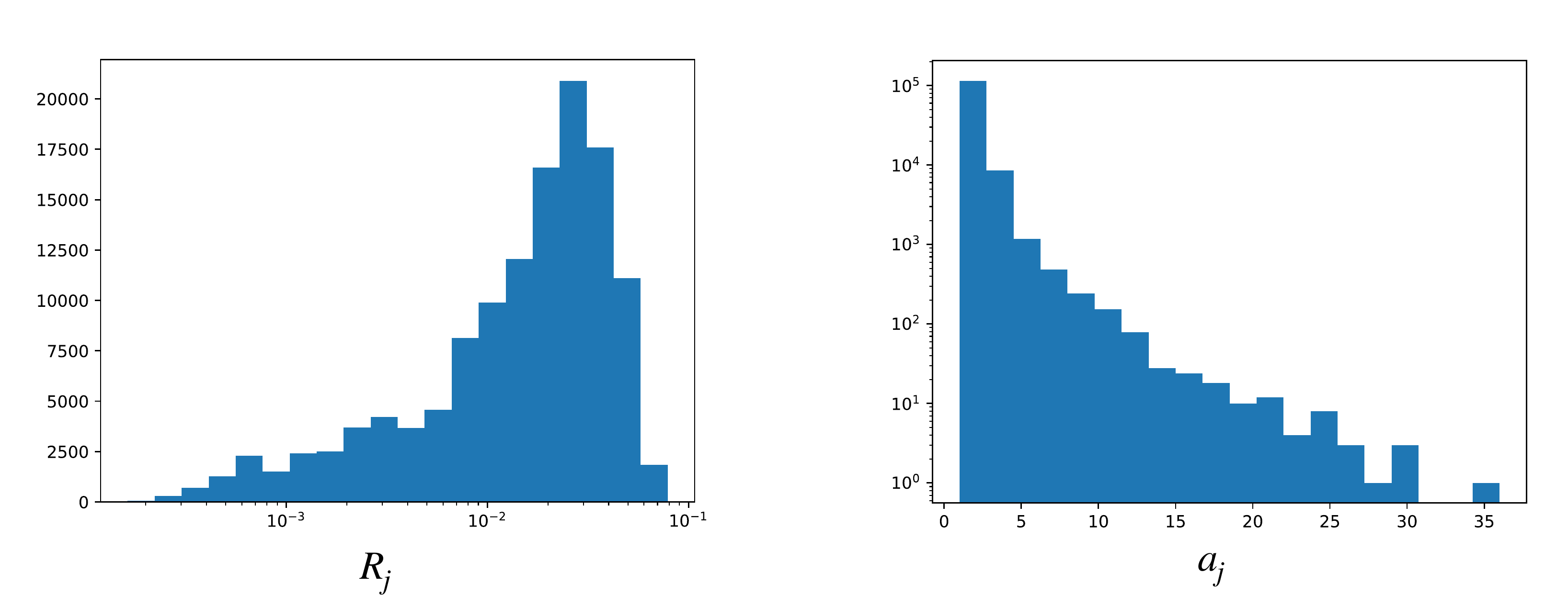}
    \end{center}
\caption{Histogram of size of patches $R_{j}$ (left), and aspect ratio
  $a_{j}$ right.}
\label{fig:hist-plane}
\end{figure}

\begin{figure}[t!]
\begin{center}
    \includegraphics[width=\linewidth]{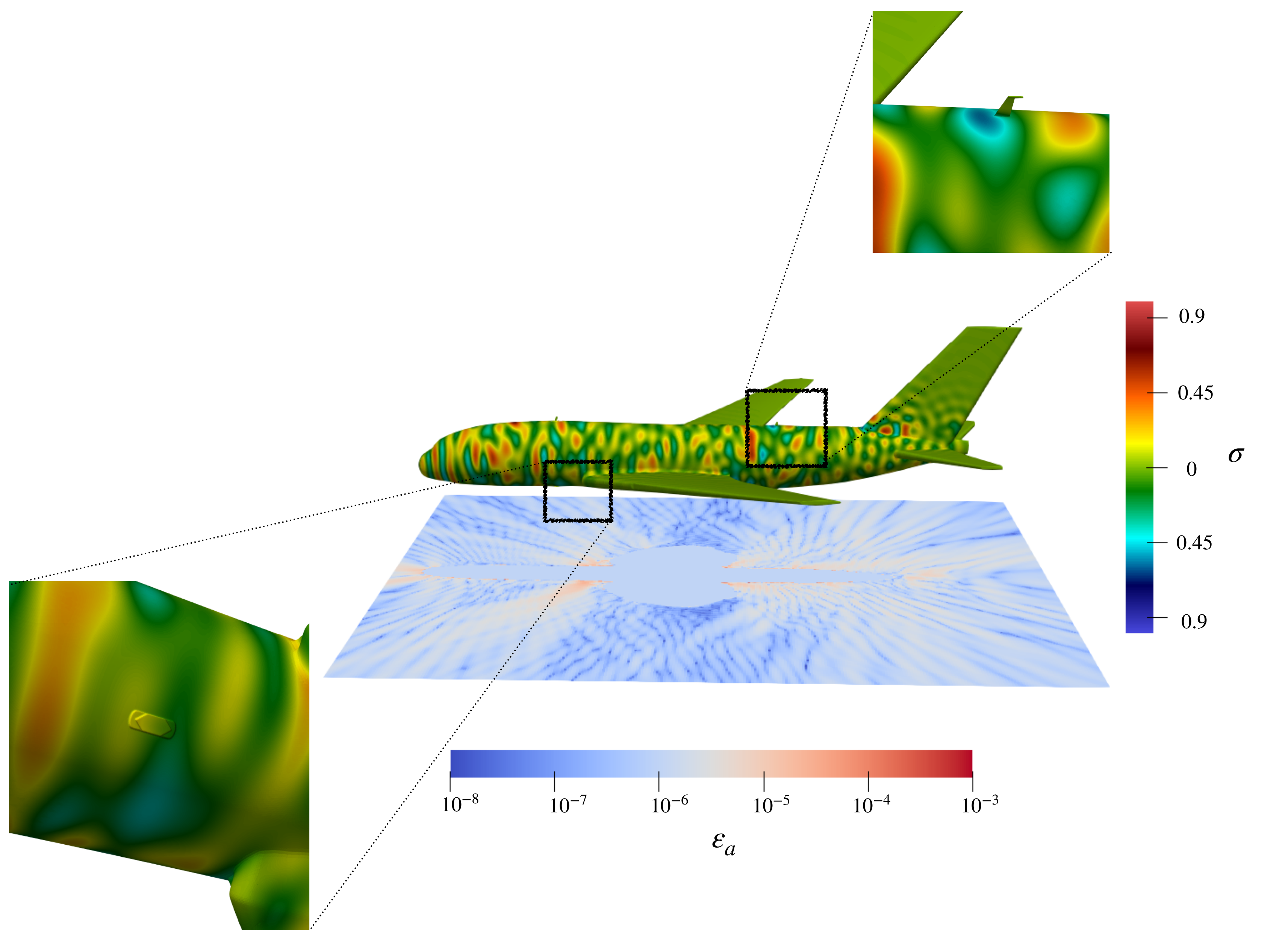}
    \end{center}
\caption{The solution $\sigma$ to~\cref{eq:plane-inteq}, and the
  relative error in the solution at a grid of $301\times301$ targets
  on a horizontal slice intersecting the wing edge.  Zoomed in views
  of one antenna (top right) and the control unit (bottom left)
  indicate the extent of the fine multiscale features.}
\label{fig:plane}
\end{figure}

In this section we demonstrate the performance of our method on a
moderate frequency acoustic scattering problem. The model airplane is
$49.3$ wavelengths long, with a wingspan of $49.2$ wavelengths and a
vertical height of $13.7$ wavelengths, which we assume has a sound-soft boundary,
satisfying Dirichlet boundary conditions (see section \ref{sec:surface}).
The plane also has several
multiscale features: 2 antennae on the top of the fuselage, and 1
\emph{control unit} on the bottom of the fuselage.  (See
\Cref{fig:plane}.)  The plane is discretized with $\Npat=125,344$, and
$p=4$ resulting in $N=1,253,440$ discretization points.  The ratio of
the largest to the smallest patch size, measured by the enclosing
sphere radius $R_{j}$ in~\cref{eq:rjdef}, is $483.9$.
In~\cref{fig:hist-plane}, we plot a histogram of the patch sizes
$R_{j}$, and the aspect ratios of the patches on the plane. The worst
case patch has an aspect ratio of $35$, but only $322$ patches out of
the total $125,344$ have an aspect ratio of greater than $10$.

In order to have an analytic reference solution, we assume the Dirichlet
boundary data $u|_\Gamma$ for the governing exterior Helmholtz equation is 
generated using a collection of $123$ interior sources, $19$ of which 
are in the tail. 
Using a combined field representation
\[
u = \cD_{k}[\sigma] - ik \, \cS_{k}[\sigma],
\]
imposing the Dirichlet condition yields the 
combined field integral equation for the unknown density~$\sigma$:
\begin{equation}
\frac{\sigma}{2} + \cD_{k}[\sigma] - ik \cS_{k}[\sigma] = u|_\Gamma \, .
\label{eq:plane-inteq}
\end{equation}


For a quadrature tolerance of $5 \times 10^{-7}$, 
$\tinit=75.86$s, the oversampling factor $\alpha=4.51$, and the memory
cost per discretization point is $m=149.6$. GMRES converged to a relative 
residual of $5 \times 10^{-7}$ in $59$ iterations, and the solution was 
obtained in $4,694$s. 
We plot the solution on a $301 \times 301$ lattice of targets on a slice 
which cuts through the wing edge, whose normal is given 
by $(0,0,1)$. For targets in the interior of the airplane,
we set the error to $1\times 10^{-6}$ since the computed solution there is 
not meaningful. The layer potential evaluation at all targets required only
$\tinit + \tquad = 104.37$s. In~\cref{fig:plane}, we plot the density 
$\sigma$ on the airplane surface and the relative error on the slice with
$301 \times 301$ targets. The maximum relative error at all targets 
is $5 \times 10^{-4}$.

\subsubsection*{Scattering from an airplane (sound-hard)}
  \label{subsec:plane2}

In this section we solve the Helmholtz equation in the exterior
of the plane, assuming Neumann boundary conditions instead.
These arise in the modeling of sound-hard scatterers \cite{colton_kress}.
We will use the following regularized combined field
integral representation of~\cite{bruno_2012} for the solution:
\begin{equation}
  \label{eq:plane-rep_hard}
    u=\cS_{k}[\sigma]+i\alpha \cD_{k}\left[ \cS_{i|k|}[\sigma] \right].
\end{equation}
Applying the boundary condition $\partial u/\partial  n=g$
along~$\Gamma$ leads to the second-kind integral equation:
\begin{equation}
  \label{eq:plane-int_hard_far}
  -\frac{\sigma}{2} + \cS'_{k}[\sigma] + i\alpha
  \cD'_{k}\left[ \cS_{i|k|}[\sigma] \right] = g.
\end{equation}
Using a well-known Calder\'on identity for the operator $\cD'_{i|k|}
\cS_{i|k|}$~\cite{nedelec}, this equation can be re-written so as to
avoid the application of the hypersingular operator~$\cD_k'$:
\begin{equation}
  \label{eq:plane-int_hard_self}
  -\left( \frac{2+i\alpha}{4} \right) \sigma
  + \cS'_{k}[\sigma] + i\alpha \left( \cD'_{k}-\cD'_{i|k|} \right)
  \left[ \cS_{i|k|}[\sigma] \right] +
  i\alpha \cS'^2_{i|k|}[\sigma] = g.
\end{equation}
A number of different possibilities for the regularizing operator are
available depending on the frequency range
(see~\cite{colton_kress_inverse, vico_2014}).  Examining the above
integral equation, it is clear that a total of four separate FMM calls
and four local quadrature corrections will be needed.

\begin{figure}[t]
\begin{center}
    \includegraphics[width=\linewidth]{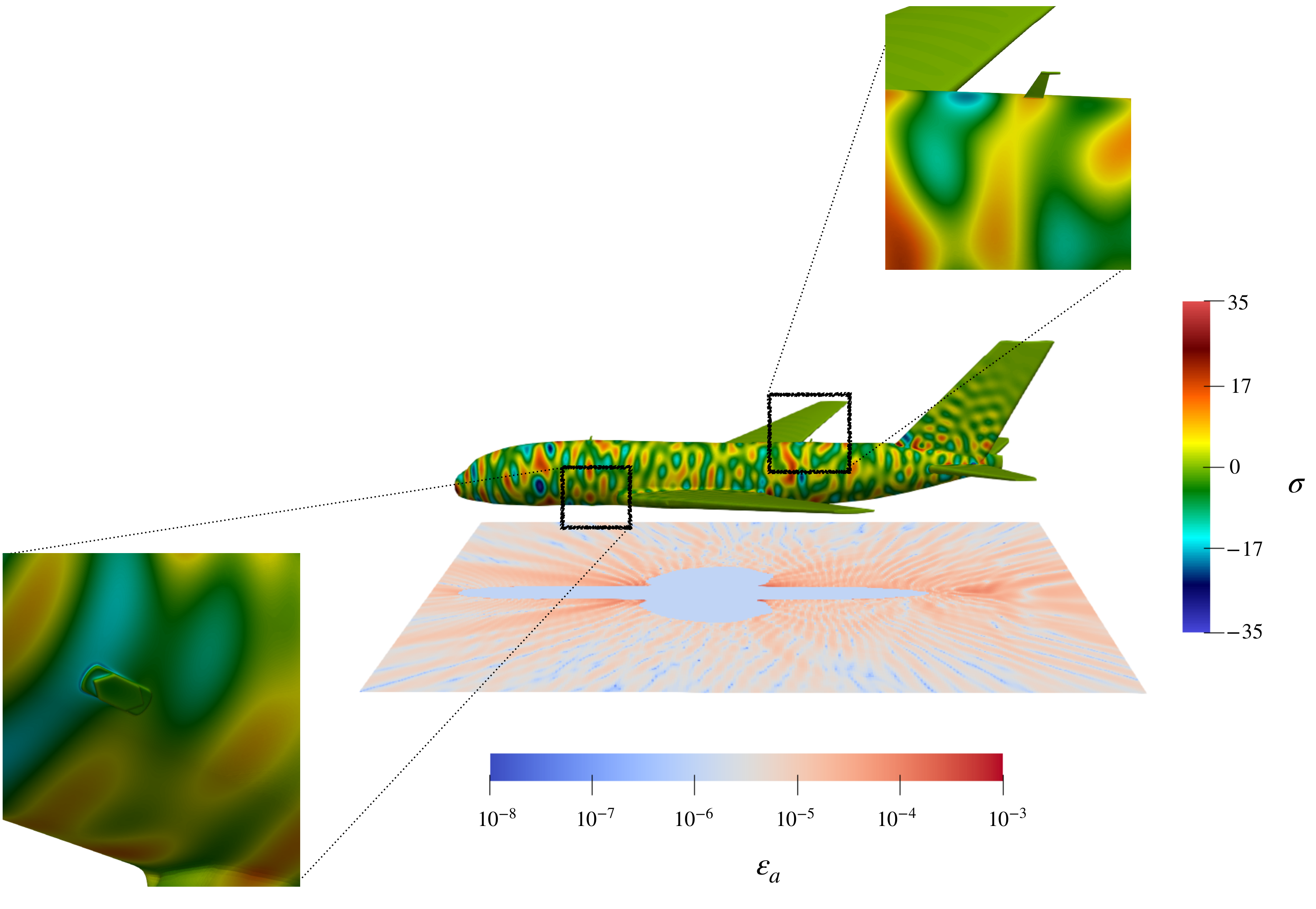}
    \end{center}
\caption{The solution~$\sigma$ to integral
  equation~\eqref{eq:plane-int_hard_self}, and the relative error in
  the solution to the PDE at a grid of~$301\times301$ targets on a horizontal
  slice intersecting the wing edge.  Zoomed in views of one antenna
  (top right) and the control unit (bottom left) indicate the extent
  of the fine multiscale features.}
\label{fig:plane-neu}
\end{figure}

In order to have an analytic reference solution, we generate the
Neumann boundary data $u|_\Gamma$ for the governing exterior Helmholtz
equation by using the same~$123$ interior sources as in the Dirichlet
problem.  For a quadrature tolerance of $5 \times 10^{-7}$,
$\tinit=387.4$s, the oversampling factor $\alpha=4.51$, and the memory
cost per discretization point is $m=598.3$. GMRES converged to a
relative residual of $5 \times 10^{-7}$ in $35$ iterations, and the
solution was obtained in $7617$s.  We plot the solution on the $301
\times 301$ lattice of targets used for the Dirichlet problem. As
before, for targets in the interior of the airplane, we set the error
to $1\times 10^{-6}$ since the computed solution there is not
meaningful. The layer potential evaluation at all targets required
only $\tinit + \tquad = 159.6$s. In~\cref{fig:plane-neu}, we plot the
induced density~$\sigma$ on the airplane surface and the relative
error on the slice with $301 \times 301$ targets. The maximum relative
error at all targets is $5 \times 10^{-4}$.

\subsubsection*{Scattering through a Fresnel lens
  (multiple sound speeds)}
\label{subsec:lens}

\begin{figure}[t]
\centering
\includegraphics[width=0.9\linewidth]{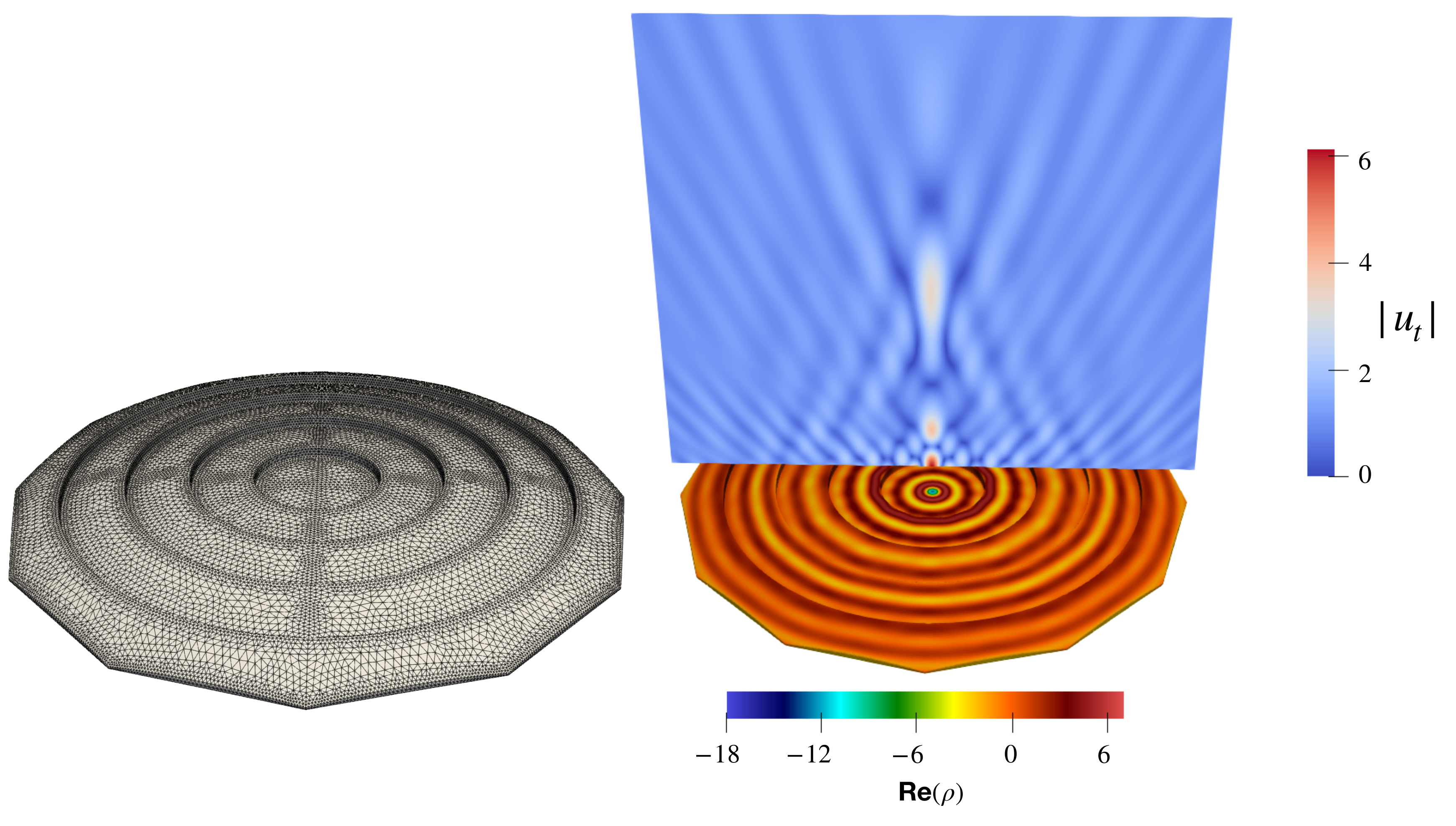}
\caption{On the left, we illustrate the triangulation of the lens surface. On the right is a plot of the real part of the solution $\rho$ to~\cref{eq:trans-sys} on the surface of the lens and the absolute value of the total field above the lens
where the central focusing of the beam is clearly visible.}
\label{fig:fresnel}
\end{figure}

In this section we solve the Helmholtz transmission problem,
i.e. scattering through media with various sound speeds, in a Fresnel
lens geometry (see~\Cref{fig:fresnel}).  The Helmholtz
parameter for the interior region (the lens)
is~$k=\omega\sqrt{\epsilon\mu}$  and for the exterior region (free
space) is $k_0=\omega\sqrt{\epsilon_0\mu_0}$.
We assume that the known incoming field~$u^{\In}$ exists only in the
exterior, so that the total field in the exterior is given by~$u_t =
u_0 + u^{\In}$ where~$u_0$ is the scattered field. In the interior,
the total field is merely $u_t = u$, with~$u$ the scattered field.

In the piecewise constant sound speed setup, we enforce the following
transmission conditions across interfaces:
\begin{equation}
  \begin{aligned}
    u_0-u&=-u^{\In}|_\Gamma, \\ \frac{1}{\epsilon_0}\frac{\partial
      u_0}{\partial n}-\frac{1}{\epsilon_1}\frac{\partial u}{\partial
      n}&=-\frac{1}{\epsilon_0}\frac{\partial u^{\In}}{\partial
      n}\Big|_\Gamma.
  \end{aligned}
\end{equation}
The scattered field in the exterior region,~$u_0$, and in the interior
region,~$u$, are represented as, respectively~\cite{colton_kress}:
\begin{equation}
  \begin{aligned}
    u_0&=\epsilon_0\cD_{k_0}[\rho]+ \epsilon_0^2\cS_{k_0}[\sigma], \\
    u&=\epsilon\cD_{k}[\rho]+ \epsilon^2\cS_{k}[\sigma].
    \end{aligned}
\end{equation}
This leads to the system of boundary integral equations along~$\Gamma$
\begin{equation}
    \begin{aligned}
    \Big(\frac{\epsilon_0+\epsilon}{2}\Big)\rho+\Big(\epsilon_0\cD_{k_0}-\epsilon\cD_{k}\Big)[\rho]+\Big(\epsilon_0^2\cS_{k_0}-\epsilon^2\cS_{k}\Big)[\sigma]&
    = -u^{\In},\\
    \Big(\frac{\epsilon_0+\epsilon}{2}\Big)\sigma-\Big(\cD'_{k_0}-\cD'_{k}\Big)[\rho]-\Big(\epsilon_0\cS'_{k_0}-\epsilon\cS'_{k}\Big)[\sigma]&=
    - \frac{1}{\epsilon_0} \frac{\partial u^{\In}}{\partial n}.
    \end{aligned}
    \label{eq:trans-sys}
\end{equation}

In the following example, the Fresnel lens has~$\epsilon=2$,~$\mu=1$ and
as usual, in free-space~$\epsilon_0=\mu_0=1$. The angular frequency
is set to be~$\omega=1+\sqrt{2}$, and the annular step size in the
lens equals 1.  Relative to the exterior wavenumber, the Fresnel lens is 19.11 wavelengths in diameter and has a height of 0.84 wavelengths.
The geometry was designed in
GiD~\cite{gid}, and a 4th-order curvilinear
mesh was constructed using the method described
in~\cite{vico2020mesh}. The mesh consists of 62,792 curvilinear
triangles, each discretized to 5th order yielding a total of 941,880
discretization points. See Figure~\ref{fig:fresnel}.

We solve the transmission problem in response to an incoming plane wave $u^{inc}=e^{ik_{0}z}$. For a quadrature tolerance of $5 \times 10^{-7}$,
$\tinit=238.2$s, the oversampling factor $\alpha=3.07$, and the memory cost per discretization point is $m=762.0$. GMRES converged to a relative residual of $5 \times 10^{-7}$ in $294$ iterations, and the
solution was obtained in $23029$s.  In~\cref{fig:fresnel} we plot the absolute value of the total field $|u_{t}| = |u_{0} + u^{inc}|$ on a $1000
\times 1000$ lattice of targets in the $yz$ plane in the exterior region.
The layer potential evaluation at all targets required only $\tinit+\tquad=81.1$s. We also plot the real part of the density $\rho$ on the surface of the lens. 

The accuracy of the solution is estimated by solving a transmission problem whose boundary data is computed using known solutions to the Helmholtz equation in the interior and exterior (the interior Helmholtz solution is the potential due to a point source in the exterior and the exterior Helmholtz solution is the potential due to  a point source in the interior). The maximum relative error for the computed solution at the same target grid as above is $8.4\times 10^{-5}$.


\section{Conclusions}
\label{sec:conclusions}

In this paper, we have presented a robust, high-order accurate method
for the evaluation of layer potentials on surfaces in three
dimensions. While our examples have focused on the Laplace and
Helmholtz equations, the underlying methodology extends naturally to
other problems in mathematical physics where the governing Green's
function is singular -- thereby requiring specialized quadrature
schemes -- but compatible with fast multipole acceleration in the far
field.

To determine the highest performance scheme, we implemented
generalized Gaussian quadrature
\cite{bremer_2012c,bremer_2013,bremer-2015,bremer}, QBX
\cite{klockner_2013,Siegel2018ALT,Wala2018,Wala2020}, and coordinate
transformation schemes similar to
\cite{bruno2001fast,malhotra19,ying}. After various code optimizations
(at least for surfaces defined as collections of curved, triangular
patches), we found that generalized Gaussian quadrature with careful
reuse of precomputed, hierarchical, adaptive interpolation tables was
most efficient, as illustrated in the preceding section. It may be
that an even better local quadrature scheme emerges in the future (in
particular, extensions of the idea presented in~\cite{wu2020corrected} look quite
promising).  In that case, as discussed in \Cref{sec:fmmcoupling},
coupling to the FMM will be fundamentally unchanged.

A useful feature of locally corrected quadrature rules (of any 
kind) is that the procedure is trivial to parallelize by 
assigning a different patch to each computational thread.
Thus, acceleration on multi-core or high performance platforms
is straightforward. 

Finally, although we limited ourselves here to evaluating layer
potentials and solving integral equations iteratively, we note that
our quadrature generation scheme and oct-tree data structure 
are compatible for use with fast direct solvers
\cite{bremer-2015,greengard-2009,ho-2012,martinsson-2005,borm_2015,
guo_2016,liu_2016,coulier_2015,minden_2016}. These require
access to small blocks of the system matrix in an effort to 
find a compressed representation  of the inverse. 

There are still a number of open questions that remain to be addressed,
including the development of rules for surfaces with edges and corners
and the coupling of layer potential codes with volume integral codes
to solve inhomogeneous or variable-coefficient partial differential
equations using integral equation methods. These are all ongoing areas of
research.

\section*{Acknowledgments}
\label{sec:ack}

We would like to thank Alex Barnett and Lise-Marie Imbert-G\'erard
for many useful discussions, and Jim Bremer and
Zydrunas Gimbutas for sharing several useful quadrature codes.
We also gratefully acknowledge the support of the NVIDIA Corporation 
with the donation of a Quadro P6000, used for some of the 
visualizations presented in this research. L. Greengard was supported in part by 
the Office of Naval Research under award number~\#N00014-18-1-2307.
M. O'Neil was supported in part by the Office of Naval Research under award
numbers~\#N00014-17-1-2059,~\#N00014-17-1-2451, and~\#N00014-18-1-2307.
F. Vico was supported in part by the Office of Naval Research under award number~\#N00014-18-2307, 
the Generalitat Valenciana under award number AICO/2019/018, and by the Spanish Ministry of Science and Innovation
(Ministerio Ciencia e Innovaci\'{o}n) under award number PID2019-107885GB-C32.
We would also like to thank the anonymous referees for many helpful comments that led to a much-improved manuscript.

 \bibliographystyle{elsarticle-num}
 \bibliography{master}
\end{document}